\documentclass[12pt]{amsart}
\usepackage{latexsym}
\usepackage{a4}
\usepackage{amssymb}
\newtheorem{thm}{Theorem}[section] 
\newtheorem{pro}[thm]{Proposition} 
\newtheorem{lem}[thm]{Lemma} 
 
\newtheorem{cor}[thm]{Corollary} 
\theoremstyle{definition}
 
\newtheorem{prb}[thm]{Problem}
\newtheorem{de}[thm]{Definition}
\numberwithin{equation}{section} 

\newcommand{\ab}[1]{{\mathbf{#1}}} 
\newcommand{\ob}[1]{{\mathbb{#1}}} 

\newcommand{\cc}[1]{\ab{#1}^c}
\newcommand{\mygamma}{{\boldsymbol \Gamma}}
\newcommand{\mydelta}{{\boldsymbol \Delta}}

\newcommand{\N}{\Bbb{ N}}

\newcommand{\setsuchthat}{\,\, \pmb{|} \,\,} 

\newcommand{\congmod}[3]{#1 \equiv #2 \,\, \left(\mbox{\rm mod } #3\right)}

\newcommand{\vb}[1]{\mathbf{#1}}

\newcommand{\LComm}[2]{\lceil #1, #2 \rceil}

\newcommand{\Clo}{\mathrm{Clo}}

\newcommand{\Comp}{\mathrm{Comp}}

\newcommand{\Con}{\mathrm{Con}}

\newcommand{\meet}{\wedge}
\newcommand{\join}{\vee}
\newcommand{\bigjoin}{\bigvee}
\newcommand{\bigmeet}{\bigwedge}
\newcommand{\Int}[2]{{I[{#1},{#2}]}}

\newcommand{\lcover}{\prec}

\newcommand{\projective}{\leftrightsquigarrow}

\newcommand{\algop}[2]{( {#1}, {#2} )}

\newcommand{\epsi}{\varepsilon}

\newcommand{\isomorphic}{\cong}
\title{Congruence lattices forcing nilpotency}
\author{Erhard Aichinger}

\address{Erhard Aichinger,
Institut f\"ur Algebra,
Johannes Kepler Universit\"at Linz,
4040 Linz,
Austria}
\email{\tt erhard@algebra.uni-linz.ac.at}


\urladdr{http://www.jku.at/algebra}
\thanks{Supported by the Austrian Science Fund (FWF):P24077.}
\date{\today}

\begin{document}
\bibliographystyle{amsalpha}
\begin{abstract}
Given a lattice $\ob{L}$ and a class $K$ of algebraic structures, we say
that $\ob{L}$ \emph{forces nilpotency} in $K$ if every algebra
$\ab{A} \in K$ whose congruence lattice $\Con (\ab{A})$ is isomorphic  
to $\ob{L}$ is nilpotent.
We describe congruence lattices that force nilpotency, supernilpotency or
solvability for some classes of algebras. For this purpose, we investigate
which commutator operations can exist on a given congruence lattice.
\end{abstract}

\maketitle
\section{Introduction} \label{sec:intro}
We look for structural properties of an algebraic structure that
are forced by the shape of its congruence lattice. In particular, we
will consider the following properties of an algebra: being \emph{abelian},
being \emph{solvable}, being \emph{nilpotent}, and being \emph{supernilpotent};
the first three of these properties were first introduced for groups, but they proved
meaningful for all algebraic structures. Examples of results in universal algebra
\cite{BS:ACIU}
involving these concepts are that in a congruence modular variety,
every abelian algebra is -- essentially -- a ring module \cite{He:AAIC, Gu:GMIC} and
that every nilpotent algebra of prime power order has a loop reduct,
permutable congruences, and generates a finitely axiomatizable variety \cite{FM:CTFC}.
For a property $p_A$  of an \emph{algebra}, we search for
a corresponding property $p_L$ of a \emph{lattice} such that
every algebra whose congruence lattice satisfies $p_L$ has the property $p_A$.
Since arbitrary algebras can be quite diverse, all our results will be applicable
only to restricted classes of algebras, such as the class
$D$ of all algebras generating congruence modular varieties.
\begin{de}
  Let $K$ be a class of universal algebras, and let $P$
  be the subclass of those algebras in $K$ that fulfil the property $p$.
  Let $\ob{L}$ be a lattice. Then $\ob{L}$ \emph{forces $p$ in $K$} if
  every algebra $\ab{A} \in K$ such that $\Con (\ab{A})$ is isomorphic to $\ob{L}$
  lies in $P$.
\end{de}
We will consider this definition first with $K := D$ and $P$
the subclass of solvable algebras in $D$. Then we could pose the following problem:
\begin{quote}
  Characterize those finite lattices that force solvability in $D$.
\end{quote}
However, among these lattices we also find those finite modular lattices
that do not appear as congruence lattices of algebras in $D$. Hence a property
$p_L$ characterizing these lattices most hold for all the ``forbidden'' finite modular lattices
that never appear as congruence lattices of an algebra in $D$. This difficulty can be avoided if we only
consider those lattices that actually are congruence lattices. To this end, for a
class $K$ of algebras, we define the class $L (K)$ by
$L (K) := \{ \ob{L} \setsuchthat \exists \ab{A} \in K : \ob{L} \cong \Con (\ab{A}) \}$
as the class of congruence lattices of algebras in $K$.
Then in the present note we will
\begin{enumerate}
\item characterize those lattices that force solvability (or supernilpotency)
    in $D$ among the lattices
  of finite height in $L (D)$;
\item characterize those lattices that force nilpotency in $G$ among the lattices
  in $L (G)$, where $G$ is the class of finite expanded groups;
 \end{enumerate}
The properties that characterize these lattices will be rather easy to state, provided that
we have some basic notions from lattice theory \cite{MMT:ALVV,Gr:GLTS} at our disposal. 
We call $\Int{\alpha}{\beta}$ a \emph{prime interval} of the lattice $\ob{L}$, write $\alpha \lcover \beta$, and
say that $\alpha$ is a \emph{subcover} of $\beta$
if $\alpha < \beta$ and the interval $\Int{\alpha}{\beta}$ is exactly the set $\{\alpha, \beta\}$.
Departing from common usage,
we call an element $\eta$ of a complete lattice \emph{meet irreducible} if $\eta < 
\bigmeet \{ \beta \setsuchthat \eta < \beta \}$, and in this
case we abbreviate 
$\bigmeet \{ \beta \setsuchthat \eta < \beta \}$ by $\eta^+$.
The set of meet irreducible elements of the complete lattice $\ob{L}$ is denoted
by $M(\ob{L})$.
For arbitrary $\alpha, \beta, \gamma, \delta \in \ob{L}$,
we write $\Int{\alpha}{\beta} \nearrow \Int{\gamma}{\delta}$ if
$\delta = \beta \join \gamma$ and $\alpha = \beta \meet{\gamma}$; 
\emph {projectivity} is the smallest equivalence on intervals 
containing $\nearrow$, and it is denoted  by $\projective$.
We first state a  description of finite lattices that force solvability.
\begin{thm} \label{thm:solv}
  Let $\ob{L}$ be a lattice of finite height that
  is the congruence lattice of some
  algebra in a congruence modular variety. Then the following are
  equivalent:
  \begin{enumerate}
  \item \label{it:s1} $\ob{L}$ forces solvability in the class
     of algebras generating congruence modular varieties.
    \item \label{it:s2}
    Every algebra $\ab{B}$ generating a congruence modular variety with
    $\Con (\ab{B}) \isomorphic \ob{L}$ is solvable.
  \item \label{it:s4} The two element lattice $\ob{B}_2$ is not a homomorphic
       image of $\ob{L}$.
  \end{enumerate}
\end{thm}
We notice that for finite algebras, the implication \eqref{it:s4}$\Rightarrow$\eqref{it:s2}
is a consequence of \cite[Theorem~7.7(2)]{HM:TSOF}.

For a prime interval
$\Int{\alpha}{\beta}$ of the  complete lattice $\ob{L}$,
 we define the element $\mygamma (\alpha, \beta)$ of $\ob{L}$ by
\[
        \mygamma (\alpha, \beta) := \bigjoin \{ \eta \in M (\ob{L}) \setsuchthat 
                           \Int{\eta}{\eta^+} \projective \Int{\alpha}{\beta} \}.
\]           
Using these elements $\mygamma (\alpha, \beta)$, we can express a condition
forcing nilpotency in finite expanded groups.
\begin{thm} \label{thm:nil}
  Let $\ob{L}$ be a lattice that is the congruence lattice of some
  finite expanded group.
  Then the following are equivalent:
  \begin{enumerate}
    \item \label{it:n1}
        $\ob{L}$ forces nilpotency in the class of finite expanded groups.
    \item \label{it:n2}
    Every finite expanded group $\ab{B}$ with
    $\Con (\ab{B}) \isomorphic \ob{L}$ is nilpotent.
  \item \label{it:n3}
    For each prime interval $\Int{\alpha}{\beta}$ of $\ob{L}$, we have
    $\mygamma (\alpha, \beta) = 1$.
  \end{enumerate}
\end{thm}

The third algebra property for which a lattice property was found is supernilpotency.
The following theorem gives a description of congruence lattices that
force supernilpotency. We say that a lattice $\ob{L}$ \emph{splits} if
it is the union of two proper subintervals, which is equivalent to
saying $\ob{L} \models \exists \, \delta, \epsi : (\delta < 1 \text{ and } \epsi > 0 \text{ and }
                                                \forall \alpha : (\alpha \le \delta \text{ or } \alpha \ge \epsi))$.
A pair $(\delta, \epsi) \in (\ob{L} \setminus \{1_A \}) \times (\ob{L} \setminus  \{0_A\})$ with
$\ob{L} = \Int{0_A}{\delta} \cup \Int{\epsi}{1_A}$ is also called a \emph{splitting pair} of $\ob{L}$.
\begin{thm} \label{thm:snp}
  Let $\ob{L}$ be a finite lattice that is the congruence lattice
  of some algebra in a congruence modular variety.
  Then the following are equivalent:
   \begin{enumerate}
     \item \label{it:su1}
        $\ob{L}$ forces supernilpotency in the class of all algebras that
        generate a congruence modular variety.
     \item \label{it:su2}
    Every algebra $\ab{B}$ in a congruence modular variety with
    $\Con (\ab{B}) \isomorphic \ob{L}$ is supernilpotent.
  \item \label{it:su3}
     $\ob{L}$ does not split.
  \end{enumerate}
\end{thm} 
The proofs of 
Theorems~\ref{thm:solv}, \ref{thm:nil}, and \ref{thm:snp} are given in Section~\ref{sec:proofs}.
Parts of these results will be proved in a purely lattice theoretic setting.
To this end, the congruence lattice of an algebra is
expanded with the binary operation of taking commutators. One obtains a new
algebraic structure called \emph{commutator lattice} which has been introduced and studied in 
\cite{Cz:AOTC,Cz:TEDC}. Section~\ref{sec:cl} contributes to the structure theory of these 
commutator lattices.

\section{Preliminaries on congruence lattices and commutators} \label{sec:prel}
When seeking to describe an algebraic structure $\ab{A} = \algop{A}{F}$,
we can find significant information in the set of its congruence relations.
These congruence relations, ordered by $\subseteq$, are a complete sublattice
of the set of equivalence relations on the set $A$; the set of congruence relations
is denoted by $\Con (\ab{A})$. For arbitrary algebras,
these  congruence relations play the role that ideals play for rings and that
normal subgroups play for groups. Commutator theory \cite{FM:CTFC} generalizes taking the
commutator subgroup of two normal subgroups to arbitrary algebraic structures
by associating
a new congruence $\gamma := [\alpha, \beta]_{\ab{A}}$ with every pair of congruences $(\alpha, \beta)$ of
$\ab{A}$.
Generalizations of the group commutator can be found, e.g., in \cite{Hi:GWMO} 
and \cite{Sc:TSOO}, but it was the work of \cite{Sm:MV,HH:ACIM,FM:CTFC} 
that led to the following definition of the \emph{term condition commutator}, which
generalizes at the same time taking the commutator subgroup $[A,B]$ of
two normal subgroups of a group, and forming the ideal product $A\cdot B$ of two
ideals of a ring.
\begin{de}[cf. {\cite[Definition~4.150]{MMT:ALVV}}]
   Let $\ab{A}$ be an algebraic structure, and let $\alpha, \beta$ be congruences of $\ab{A}$.
   Then the commmutator $\gamma := [\alpha, \beta]_{\ab{A}}$ is defined as the 
   intersection of all congruence relations $\delta$ of $\ab{A}$ such that
   for all $n \in \N$, for all $(n+1)$-ary term functions $t$ of $\ab{A}$, and
   for all $(a,b) \in \alpha$ and $(c_1, d_1),\ldots, (c_n, d_n) \in \beta$ with
   \(
      (t(a,c_1,\ldots,c_n), t(a, d_1,\ldots, d_n)) \in \delta \)
    we have
   \(
      (t(b,c_1,\ldots,c_n), t(b, d_1,\ldots, d_n)) \in \delta.
   \)
\end{de}
Defined for arbitrary
algebras, commutators have proved most useful for algebras with a modular 
congruence lattice, and hence we will restrict ourselves to
such algebras, or, in decreasing steps of generality, 
to algebras in congruence modular varieties, 
to algebras in congruence permutable
varieties, or to expanded groups. 
In congruence permutable varieties, the term condition commutator
admits the following description, which resembles the ideal product
defined in \cite{Sc:TSOO}.
\begin{lem}[cf. {\cite[Corollary~6.10]{AM:SAOH}}]
   Let $\ab{A}$ be an algebra in a congruence permutable variety, and let
   $\alpha, \beta$ be congruences of $\ab{A}$. Then the congruence $[\alpha, \beta]_{\ab{A}}$
   is generated as a congruence of $\ab{A}$ by $\{ (c (a_1, b_1), c(a_2, b_2)) \setsuchthat
    (a_1, a_2) \in \alpha,\, (b_1, b_2) \in \beta, \, 
    c \text{ is a binary polynomial function of } \ab{A}
    \text{ with } 
    c(a_1, b_1) = c(a_1, b_2) = c(a_2, b_1)
     \}$. 
\end{lem}
From the congruence lattice
and the commutator operation of a finite algebra in a congruence modular variety, 
one can, e.g., determine whether the algebra generates a residually small variety
\cite[Theorem~10.15]{FM:CTFC} or whether every homomorphic image of an algebra
in a congruence permutable variety is affine complete 
\cite[Proposition~5.2]{Ai:OHAH}. 
Starting from the
commutator operation on congruences,
it is possible to define the \emph{derived series} $(\gamma_n)_{n \in \N}$ and the
\emph{lower central series} $(\lambda_n)_{n \in \N}$ of congruences of the algebra $\ab{A}$ by
$\gamma_1 = \lambda_1 = 1_A$, and the recursion
$\gamma_{n+1} = [\gamma_n, \gamma_n]$ and $\lambda_{n+1} = [1_A, \lambda_n]$ for $n \in \N$.
An algebra in a congruence modular  variety
is called \emph{solvable} (cf. \cite[Definition~3.6(3)]{HM:TSOF}) if there is $m$ with $\gamma_m = 0_A$, and
\emph{nilpotent} (cf. \cite[p.69 before Lemma~7.3]{FM:CTFC}) if there is $k$ with $\lambda_k = 0_A$. 
A.\ Bulatov \cite{Bu:OTNO} introduced a generalization of the binary commutator operation
by associating a congruence $[\alpha_1, \ldots, \alpha_n]$ with every finite sequence
of congruences; $[\alpha_1,\ldots, \alpha_n]$ is called a \emph{higher commutator}.
In congruence modular varieties, the higher commutator operations enjoy
certain properties, such as monotonicity, symmetry, and distributivity with respect to joins; the
validity of some of these properties was established only recently in \cite{Mo:HCTF}.
If an algebra has an $m \in \N$ such that $[\alpha_1,\ldots,\alpha_n] = 0$ 
whenever $n > m$, then the algebra is
called \emph{supernilpotent}. Every supernilpotent algebra in a congruence modular
variety is nilpotent: for congruence permutable varieties, this was proved in \cite[Corollary~6.15]{AM:SAOH},
and for congruence modular varieties, it follows from properties (4) and (8) of higher commutators given
in \cite{Mo:HCTF}, which are called (HC4) and (HC8) in \cite[p.\ 860]{AM:SOCO} and in \cite{AM:SAOH}.
Supernilpotency admits the following combinatorial description: 
a finite algebra $\ab{A}$ in a congruence modular variety is supernilpotent if and only if
there exists a polynomial $p$ such that the $n$-generated
free algebra in the variety generated by $\ab{A}$ has at most $2^{p(n)}$ elements.
A self-contained version of this description
for the case that $\ab{A}$ is an expanded group has been given in 
 Section~4 of \cite{Ai:OTDD};
the general result follows from a combination of 
\cite[Theorem~9.18 and Lemma~12.4]{HM:TSOF}, \cite[Theorem 6.2, Corollary 7.5, Theorem 14.2]{FM:CTFC},
the notion of \emph{rank} from \cite[p.\ 179]{Ke:CMVW}, \cite[Lemma~7.5]{AM:SAOH},
the proof of Theorem~1 of \cite{BB:FSON}, and the generalization of the properties
of higher commutators from congruence permutable to congruence modular varieties
in \cite{Mo:HCTF}; since we will not make use of the combinatorial description
of supernilpotency in this paper, we abstain from a thorough discussion.
The definitions of binary commutators, solvability and nilpotency are compatible with
the classic use of these notions in group theory (cf. \cite[Exercise~4.156(11)]{MMT:ALVV}).  

Given the congruence lattice of an 
algebra, it is therefore interesting to know what the possible choices of the commutator operations
are.
Certain limitations are imposed by the laws
$[x, y] \approx [y, x]$, $[x,y] \le x \meet y$, $x \le y \rightarrow [x,z] \le [y,z]$,
$[x \join y, z] \approx [x, z] \join [y, z]$ that are satisfied by
every structure $\algop{\Con(\ab{A})}{\meet, \join, [.,.]_{\ab{A}}}$ arising from an
algebra $\ab{A}$ in a congruence modular variety.
It is easy to see that on the five element lattice $\ob{M}_3$, the
constant operation $[x,y]=0$ is the only such  operation definable
on this lattice; this imposes structural consequences on algebras with
such a congruence lattice \cite[Lemma~4.153]{MMT:ALVV}. Conditions on the higher commutator operations
that are imposed by the shape of the congruence lattice are given in 
\cite[Lemma~3.3]{AM:SOCO}.

Let us now  briefly review some properties of the commutator operations in
congruence modular varieties.
These properties are proved in Chapters~3 and~4 of
\cite{FM:CTFC}.
\begin{lem} \label{lem:commprop} 
Let $\ab{A}$ be an algebra in a congruence modular variety,
and let $\alpha, \alpha_1, \beta, \beta_1 \in \Con (\ab{A})$. Then
$[\alpha, \beta] = [\beta, \alpha] \le \alpha \meet \beta$, 
$[ \alpha \join \alpha_1, \beta] = [\alpha, \beta] \join [\alpha_1, \beta]$, and
if $\alpha \le \alpha_1$ and $\beta \le \beta_1$, then $[\alpha, \beta] \le [\alpha_1, \beta_1]$.
If $(\alpha_i)_{i \in I}$ is a family of congruences of $\ab{A}$, we also 
have $\bigjoin_{i \in I} [\alpha_i, \beta] = [\bigjoin_{i \in I} \alpha_i, \beta]$.
\end{lem}
The proofs of some of these properties are by no means obvious
and  require skilful manipulations with Day terms \cite{Da:ACOM, FM:CTFC}. The proofs become
easier when restricting to congruence permutable varieties, and
some of these properties have been proposed as exercises in \cite{MMT:ALVV}.
The introductory chapter of \cite{Ai:TPFO2} provides solutions to some
of these exercises, as does \cite{AM:SAOH}.
\begin{lem} \label{lem:commext}
   Let $\ab{A} = \algop{A}{F}$ be an algebra in a congruence modular variety,
   and let $\ab{B} = \algop{A}{F \cup G}$ be an expansion of $\ab{A}$.
   Then for all $\alpha, \beta \in \Con (\ab{B})$, we have
   $[\alpha, \beta]_{\ab{A}} \subseteq [\alpha, \beta]_{\ab{B}}$.
   Furthermore, if $\ab{B}$ is solvable, then $\ab{A}$ is solvable,
   and if $\ab{B}$ is nilpotent, then $\ab{A}$ is nilpotent.
\end{lem}
\emph{Proof:} 
 Using the definition of the commutator $[\alpha, \beta]_{\ab{B}}$ by the term condition, 
 we obtain that $\alpha$ centralizes $\beta$ modulo  $[\alpha, \beta]_{\ab{B}}$ in $\ab{B}$ (cf.
 \cite[Definition~4.148]{MMT:ALVV}).
 Since $\Clo (\ab{A}) \subseteq \Clo (\ab{B})$, $\alpha$ centralizes $\beta$ modulo
  $[\alpha, \beta]_{\ab{B}}$ in $\ab{A}$. Hence
  $[\alpha, \beta]_{\ab{A}} \le [\alpha, \beta]_{\ab{B}}$.
  Let $(\gamma^{\ab{A}}_n)_{n \in \N}$, $(\gamma^{\ab{B}}_n)_{n \in \N}$ be the derived series of $\ab{A}$ and $\ab{B}$, resp.
  Then for each $n \in \N$, we have $\gamma^{\ab{A}}_n \le \gamma^{\ab{B}}_n$, which is proved
  by induction using $\gamma^{\ab{A}}_{n+1} = [\gamma^{\ab{A}}_n, \gamma^{\ab{A}}_n]_{\ab{A}} \le
  [\gamma^{\ab{B}}_n, \gamma^{\ab{B}}_n]_{\ab{A}} \le [\gamma^{\ab{B}}_n, \gamma^{\ab{B}}_n]_{\ab{B}} 
  = \gamma^{\ab{B}}_{n+1}$ as the induction step. 
  Hence if $\ab{B}$ is solvable, then so is $\ab{A}$. The proof for nilpotency is similar. \qed

For an algebra $\ab{A}$ in a congruence modular variety and
$\alpha, \beta \in \Con (\ab{A})$, we define $(\alpha : \beta)_{\ab{A}}$
as the largest $\gamma \in \Con (\ab{A})$ with 
$[\gamma, \beta]_{\ab{A}} \le \alpha$. We omit the subscript when the
algebra is clear from the context.
When interpreting commutator theory
in group theory, $(\alpha: \beta)$ corresponds to the centralizer
$C_G (B/A)$, where $B$ and $A$ are the normal subgroups corresponding
to $\beta$ and $\alpha$. Therefore, we will call $(\alpha:\beta)$ 
the \emph{centralizer} of $\beta$ over $\alpha$. We note
that Proposition~4.2 of \cite{FM:CTFC} guarantees that this definition
is consistent with \cite[Definition~4.150]{MMT:ALVV}.
For all $\alpha, \beta, \gamma \in \Con (\ab{A})$, we
have $[\gamma, \beta] \le \alpha$ if and only if
$\gamma \le (\alpha:\beta)$; thus each of the operations $[.,.]$
and $(.:.)$ fully determines the other.

Often, we will not use any properties of the
binary commutator operation other than its mere definition by the term condition
\cite[Definition~4.150]{MMT:ALVV} and the  properties that are  stated in Lemma~\ref{lem:commprop}. 
Hence it is useful to see what can be derived from
these conditions alone; such an investigation was started
in \cite{Cz:AOTC}.
  
\section{Preliminaries on commutator lattices} 
 In \cite{Cz:AOTC}, J.\ Czelakowski defined \emph{commutator lattices}.
 These algebraic structures capture the properties
 of the structure  $(\Con (\ab{A}), \join, \meet, [.,.]_{\ab{A}})$ that
 is constructed by expanding the congruence lattice of an algebra $\ab{A}$ in a 
 congruence modular variety  with the binary
 operation of taking commutators.
\begin{de}[{\cite[Definition~1.1]{Cz:AOTC}}]
     An algebraic structure $\ab{L} = (\ob{L}, \join, \meet, [.,.])$ is a \emph{commutator lattice}
     if $(\ob{L}, \join, \meet)$ is a complete lattice, and 
     for all $x,y \in \ob{L}$ and for all families $(x_i)_{i \in I}$ from $\ob{L}$,
     we have
     $[x,y] = [y,x]$, $[x, y] \le x \meet y$, and
     $[\bigjoin_{i \in I} x_i, y] = \bigjoin [x_i, y]$. In this case, we
    call $[.,.]$ a \emph{commutator multiplication} on the lattice $(\ob{L}, \join, \meet)$.
\end{de}
The guiding example of this definition comes from congruences and commutators. 
In fact,
we may restate Lemma~\ref{lem:commprop} as follows:
\begin{pro} \label{pro:coniscl}
    Let $\ab{A}$ be an algebra that generates a congruence modular variety,
    let $(\ob{L}, \join, \meet) := (\Con (\ab{A}), \join, \cap)$ be the congruence
    lattice of $\ab{A}$, and for $\alpha, \beta \in \Con (\ab{A})$, let
    $[\alpha, \beta]_{\ab{A}}$ denote the term condition commutator 
    of $\alpha$ and $\beta$ as defined in \cite[Definition~4.150]{MMT:ALVV}.
    Then $(\ob{L}, \join, \meet, [.,.]_{\ab{A}})$ is a commutator lattice.
\end{pro}
\emph{Proof:}
    \cite[Proposition~4.3]{FM:CTFC}. \qed

 It is a consequence of the distributivity of $[.,.]$ with respect to joins that
 the operation $[.,.]$ is monotonic with respect to the order of the lattice.
  An important operation that comes along with a commutator lattice $\ab{L}$ is that of \emph{residuation}.  
  For $x, y \in \ob{L}$, we define 
  \begin{equation} \label{eq:res}
      (x:y) := \bigjoin \{ z \in \ob{L} \setsuchthat [z,y] \le x \} 
  \end{equation}
  and call $(.:.)$ the \emph{residuation operation} associated with $\ab{L}$.
  We notice that in \cite{Cz:AOTC}, $(x:y)$ is denoted by $y \rightarrow x$;
  our notation comes from the interpretation of $(x:y)$ as the
  centralizer of $y$ over $x$ in \cite{FM:CTFC}.
  In the following lemma, we state some properties of the residuation operation.
\begin{lem} \label{lem:res}
     Let $\ab{L} = (\ob{L}, \join, \meet, [.,.])$ be a commutator lattice,
     and let $(.:.)$ be the residuation operation associated with $\ab{L}$.
      Then for all $x, y, z \in \ob{L}$ and for all families $(x_i)_{i \in I}$ from $\ob{L}$,
      we have: 
     \begin{enumerate}
         \item \label{it:l0} $[z, y] \le x$ if and only if $z \le (x : y)$,
         \item \label{it:l6} $[(x:y),\, y] \le x$,
         \item \label{it:l1} $(\bigmeet_{i \in I} x_i : y) = \bigmeet_{i \in I} (x_i : y)$,
         \item \label{it:l2} $(x : \bigjoin_{i \in I} y_i) = \bigmeet_{i \in I} (x : y_i)$,
         \item \label{it:l3} $(x : y) \ge x$, 
         \item \label{it:l4} $(x : x) = 1$, 
         \item \label{it:l5} $(x : (x : y)) \ge y$.
      \end{enumerate}
\end{lem}
   \emph{Proof:} 
     \eqref{it:l0} The ``only if''-direction is an immediate consequence
       of the definition of $(x:y)$. For ``if''-direction, we assume
       $z \le (x:y)$ and compute
       $[z, y] \le [(x:y), y] = [ 
            \bigjoin \{ z_1 \in \ob{L}  \setsuchthat [z_1, y] \le x \},
            y
                                ]
            = \bigjoin \{ [z_1, y] \setsuchthat z_1 \in \ob{L}, [z_1, y] \le x \}
            \le x$.

      \eqref{it:l6} is a consequence of~\eqref{it:l0}.

            \eqref{it:l1} For $\le$, we let $j \in I$ and notice that 
            using~\eqref{it:l6}, we have $[(\bigmeet_{i \in I} x_i : y), y] \le
             \bigmeet_{i \in I} x_i \le x_j$, and therefore
             $(\bigmeet_{i \in I} x_i : y) \le (x_j : y)$. 
             Hence $(\bigmeet_{i \in I} x_i : y) \le \bigmeet_{i \in I} (x_i : y)$.
             For $\ge$, we let $j \in J$ and compute
             $[\bigmeet_{i \in I} (x_i : y), y] \le
              [ (x_j : y), y ] \le x_j$. Hence 
              $[\bigmeet_{i \in I} (x_i : y), y] \le \bigmeet_{i \in I} x_i$,
              which implies $\bigmeet_{i \in I} (x_i : y) \le (\bigmeet_{i \in I} x_i : y)$. 

      \eqref{it:l2} For $\le$, we let $j \in I$ and compute
               $[(x : \bigjoin_{i \in I} y_i), y_j] \le
                [(x : \bigjoin_{i \in I} y_i), \bigjoin_{i \in I} y_i] \le x$,
               which implies
               $(x : \bigjoin_{i \in I} y_i) \le (x : y_j)$, and therefore
               $(x : \bigjoin_{i \in I} y_i) \le \bigmeet_{i \in I} (x : y_i)$.
               For $\ge$, we compute
               $[\bigmeet_{i \in I} (x : y_i), \bigjoin_{k \in I} y_k] =
                \bigjoin_{k \in I} [\bigmeet_{i \in I} (x : y_i), y_k] \le
                \bigjoin_{k \in I} [(x : y_k), y_k] \le x$, which implies
                 $\bigmeet_{i \in I} (x : y_i) \le  (x : \bigjoin_{k \in I} y_k)$. 
          
        \eqref{it:l3} Since $[x, y] \le x$, we have $x \le (x : y)$.

        \eqref{it:l4} Since $[1, x] \le x$, we have $1 \le (x : x)$.

        \eqref{it:l5} We have $[y, (x:y)] = [(x:y), y] \le x$, and therefore
                      $y \le (x: (x:y))$. \qed                                     

    In fact, the properties~\eqref{it:l1}-\eqref{it:l5} in Lemma~\ref{lem:res} are
    equivalent to the properties (a)-(e) listed in \cite[Theorem~3.1]{Cz:AOTC},
    and therefore provide a different axiomatization of possible residuation operations of
    commutator lattices. The following lemma is an abstraction of \cite[Chapter~9, Exercise~4]{FM:CTFC}.
\begin{lem} \label{lem:centproj}
  Let $\algop{\ob{L}}{\join, \meet, [.,.]}$ be a commutator lattice,
  and let $(.:.)$ its associated residuation.
  Let $\alpha, \beta, \gamma, \delta \in \ob{L}$ 
such that $\alpha \le \beta$,
$\gamma \le \delta$, and
$\Int{\alpha}{\beta} \projective \Int{\gamma}{\delta}$.
Then
\begin{enumerate}
   \item $(\alpha : \beta) = (\gamma : \delta)$, and
   \item $[\beta , \beta] \le \alpha$ if and only if $[\delta, \delta] \le \gamma$.
\end{enumerate}
\end{lem}
\emph{Proof:}
We assume $\Int{\alpha}{\beta} \nearrow \Int{\gamma}{\delta}$.
Then using Lemma~\ref{lem:res}, we obtain
   $(\gamma : \delta) = 
    (\gamma : \beta \join \gamma) =
    (\gamma : \beta) \meet (\gamma : \gamma) =
    (\gamma : \beta) \meet 1 =
    (\gamma : \beta) \meet (\beta : \beta) =
    (\gamma \meet \beta : \beta) =
    (\alpha : \beta)$.                         
For the second item, we first assume that
$[\beta, \beta] \le \alpha$. Then $[\delta, \delta] = [\beta \join \gamma, \beta \join \gamma]
 = [\beta, \beta] \join [\beta, \gamma] \join [\gamma, \gamma] \le \alpha \join \gamma \join \gamma
 = \gamma$. Conversely, if $[\delta, \delta] \le \gamma$, then $[\beta, \beta] \le \gamma$, and since
$[\beta,\beta] \le \beta$, we obtain $[\beta, \beta] \le \gamma \meet \beta = \alpha$. \qed

Let
$\alpha, \beta \in \ob{L}$ with $\alpha \lcover \beta$.
The next lemma states that $\mygamma (\alpha, \beta) =
\bigjoin \{ \eta \in M (\ob{L}) \setsuchthat 
                           \Int{\eta}{\eta^+} \projective \Int{\alpha}{\beta} \}$
is a lower bound for the
residuum $(\alpha : \beta)$.
\begin{lem} \label{lem:lb}
  Let $\ab{L} = \algop{\ob{L}}{\join, \meet, [.,.]}$ be a commutator lattice,
  and let $(.:.)$ be its associated residuation.
  Let $\alpha, \beta \in \ob{L}$ be such that $\alpha \lcover \beta$,
  and let
  $\mygamma (\alpha, \beta) =
\bigjoin \{ \eta \in M (\ob{L}) \setsuchthat 
                           \Int{\eta}{\eta^+} \projective \Int{\alpha}{\beta} \}$.
   Then $\mygamma (\alpha, \beta) \le (\alpha : \beta)$.
\end{lem}
\emph{Proof:}
 For every $\eta \in M (\ob{L})$ with
 $\Int{\eta}{\eta^+} \projective \Int{\alpha}{\beta}$, Lemma~\ref{lem:res}\eqref{it:l3}
 and Lemma~\ref{lem:centproj} yield $\eta \le (\eta : \eta^+) = (\alpha : \beta)$.
 Therefore $\mygamma (\alpha, \beta) \le (\alpha : \beta)$. \qed

\section{Tools from lattice theory}
In constructing commutator multiplications on given lattices, we will need 
some techniques from lattice theory. We will often work in \emph{algebraic}
lattices \cite[Definition~2.15]{MMT:ALVV}, and we call 
a lattice \emph{bialgebraic} if the lattice and its dual are both algebraic; 
for example,
every lattice of finite height is bialgebraic.
For our purpose, the most important fact in algebraic lattices
is that every element is the meet of meet irreducible elements \cite[Theorem~2.19]{MMT:ALVV}.
 For any complete lattice $\ob{L}$,
$M(\ob{L})$ denotes the set of 
meet irreducible elements of $\ob{L}$, and by
$J(\ob{L}) := \{ \rho \in \ob{L} \setsuchthat \rho > \bigjoin \{ \alpha \in \ob{L} \setsuchthat \alpha < \rho \} \}$,
we denote
the set of join irreducible elements
of $\ob{L}$. The unique subcover
of a  join irreducible element $\beta$ is denoted
by $\beta^-$, and $\beta$ is called a \emph{lonesome join irreducible
element} of $\ob{L}$ if
$\{\rho \in J(\ob{L}) \setsuchthat \Int{\rho^-}{\rho} \projective
 \Int{\beta^-}{\beta}\} = \{\beta\}$; a meet irreducible element $\eta$ of
$\ob{L}$ is called a \emph{lonesome meet irreducible element} if
 $\{\varphi \in M(\ob{L}) \setsuchthat \Int{\varphi}{\varphi^+} \projective
                                       \Int{\eta}{\eta^+}\} = \{\eta\}$.
For $\alpha, \beta, \gamma, \delta \in \ob{L}$ with $\alpha \le \beta$ and $\gamma \le \delta$,
 we say that
$\Int{\alpha}{\beta}$ \emph{projects into} $\Int{\gamma}{\delta}$ if there
are $\alpha_1, \beta_1 \in \ob{L}$ with $\gamma \le \alpha_1 \le \beta_1 \le \delta$
such that $\Int{\alpha}{\beta} \projective \Int{\alpha_1}{\beta_1}$.

The following proposition collects some well known facts on projectivity.
 \begin{pro} \label{pro:projprops}
    Let $\ob{L}$ be a complete lattice.
    \begin{enumerate}
       \item \label{it:mi1} If $\ob{L}$ is algebraic and  $\alpha,  \beta \in \ob{L}$ are such that $\beta \not\le \alpha$, there
             is $\eta \in M(\ob{L})$ such that $\alpha \le \eta$ and $\beta \not\le \eta$. 
       \item \label{it:mi2} If $\ob{L}$ is modular and $\alpha, \beta \in \ob{L}$ are such that $\alpha \lcover \beta$ and $\eta \in M(\ob{L})$
             satisfies $\alpha \le \eta$ and $\beta \not\le \eta$, then $\Int{\alpha}{\beta} \nearrow
             \Int{\eta}{\eta^+}$.
       \item \label{it:mi2a} If $\ob{L}$ is modular, and $\beta \in \ob{L}$ and $\eta \in M (\ob{L})$ are such that
             $\eta \not\ge \beta$, then $\Int{\eta}{\eta^+} \searrow \Int{\eta \meet \beta}{\eta^+ \meet \beta}$.
             Dually, if $\ob{L}$ is modular, $\gamma \in \ob{L}$ and $\rho \in J(\ob{L})$ are such that
             $\rho \not\le \gamma$, then $\Int{\rho^-}{\rho} \nearrow \Int{\rho^- \join \gamma}{\rho \join \gamma}$.
       \item \label{it:mi3} If $\ob{L}$ is algebraic and  modular, and $\beta \in J(\ob{L})$ and $\gamma \in \ob{L}$
             are such that $\beta \not\le \gamma$, then there exists $\eta \in M(\ob{L})$ are such that 
             $\gamma \le \eta$ and $\Int{\beta^-}{\beta} \nearrow \Int{\eta}{\eta^+}$.
       \item \label{it:mi4} If $\ob{L}$ is modular and $a,b,x,y \in \ob{L}$ are such that 
             $x \meet y \le a \lcover b \le x \join y$, then $\Int{a}{b}$ projects into 
             $\Int{x}{x \join y}$ or into $\Int{y}{x \join y}$.
     \end{enumerate}
\end{pro}
\emph{Proof:}
    \eqref{it:mi1} By \cite[Theorem~2.19]{MMT:ALVV}, $\alpha = \bigmeet \{\psi \in M (\ob{L}) \setsuchthat
                   \psi \ge \alpha \}$. If~\eqref{it:mi1} fails, then 
                   $\{\psi \in M (\ob{L}) \setsuchthat \psi \ge \alpha \} \subseteq 
                    \{\psi \in M (\ob{L}) \setsuchthat \psi \ge \beta \}$,
                   and thus $\alpha =  \bigmeet \{\psi \in M (\ob{L}) \setsuchthat \psi \ge \alpha \} \ge
                    \bigmeet \{\psi \in M (\ob{L}) \setsuchthat \psi \ge \beta \} \ge \beta$, contradicting 
                    the assumptions.
 
    \eqref{it:mi2} Since $\beta > \eta \meet \beta \ge \alpha$, we have $\eta \meet \beta = \alpha$,
                    and from $\eta \join \beta > \eta$, we obtain $\eta \join \beta \ge \eta^+$.
                    Now suppose $\eta \join \beta \not\le \eta^+$. Then $\beta \not\le \eta^+$,
                    and therefore $\beta > \eta^+ \meet \beta \ge \alpha$. Hence
                    $\eta^+ \meet \beta = \alpha$, which implies 
                    $\eta = \eta \join (\beta \meet \eta^+)$. By modularity, we have
                    $\eta \join (\beta \meet \eta^+) = (\eta \join \beta) \meet \eta^+ = \eta^+$.
                    The contradiction $\eta = \eta^+$ completes the proof of $\eta \join \beta \le \eta^+$,
                    and therefore $\eta \join \beta = \eta^+$.

     \eqref{it:mi2a} By modularity, we have $\eta \join (\beta \meet \eta^+) = (\eta \join \beta) \meet \eta^+ =
                     \eta^+$, which proves $\Int{\beta \meet \eta}{\beta \meet \eta^+} \nearrow
                                            \Int{\eta}{\eta^+}$. The statement on join irreducible
                    elements follows from a dual argument.
   
   \eqref{it:mi3}   We first show 
   \begin{equation} \label{eq:sp1}
       \beta \not\le \gamma \join \beta^-.
   \end{equation}
    Suppose $\beta \le \gamma \join \beta^-$. Then by modularity,
    $\beta = \beta \meet (\gamma \join \beta^-) =
            (\beta \meet \gamma) \join \beta^-$. Since $\beta \not\le \gamma$,
    we have
    $\beta \meet \gamma \le \beta^-$ and thus  $(\beta \meet \gamma) \join \beta^- =
     \beta^-$. The contradiction $\beta = \beta^-$ establishes~\eqref{eq:sp1}.
    Using~\eqref{eq:sp1} and item~\eqref{it:mi1}, we find $\eta \in M(\ob{L})$ such that
    $\gamma \join \beta^- \le \eta$ and $\beta \not\le \eta$.
    Since $\beta^- \le \eta$, item~\eqref{it:mi2} yields $\Int{\beta^-}{\beta} \nearrow
    \Int{\eta}{\eta^+}$. 

       \eqref{it:mi4} If $b \join x = a \join x$ and $b \meet x = a \meet x$, then
                      $a = a \join (x \meet a) = a \join (x \meet b) = (a \join x) \meet b
                         = (b \join x) \meet b = b$.
                      Hence $b \join x > a \join x$ or $b \meet x > a \meet x$. \
                      In the case $b \join x > a \join x$, we first observe that
                      then $a \join x \not\ge b$. Hence
                      $a \le (a \join x) \meet b < b$, and therefore
                      $(a \join x) \meet b = a$, which implies $\Int{a}{b} \nearrow
                      \Int{a \join x}{b \join x}$, and therefore $\Int{a}{b}$ projects
                      into $\Int{x}{x \join y}$.
                      In the case $b \meet x > a \meet x$, we have 
                      $b \meet x \not\le a$, and therefore
                      $a < a \join (x \meet b) \le b$, which implies $a \join (x \meet b) = b$,
                      and therefore $\Int{a \meet x}{b \meet x} \nearrow \Int{a}{b}$.
                      Thus $\Int{a}{b}$ projects into $\Int{x \meet y}{x}$, and therefore
                      by Dedekind's transposition principle \cite[2.27]{MMT:ALVV} into
                      $\Int{y}{x \meet y}$. \qed

Projectivity plays an important role in the description of congruence generation in 
lattices. In a complete lattice $\ob{L}$, a relation $\Phi$ on $\ob{L}$ 
is called a \emph{complete congruence} of $\ob{L}$
if $\Phi$ is an equivalence relation on $\ob{L}$, and for all families
$(x_i)_{i \in I}$ and $(y_i)_{i \in I}$ from $\ob{L}$, we have
$(\forall i \in I : (x_i, y_i) \in \Phi) \Rightarrow ((\bigjoin_{i \in I} x_i, \bigjoin_{i \in I} y_i) \in \Phi
                                                     \text{ and }
                                                     (\bigmeet_{i \in I} x_i, \bigmeet_{i \in I} y_i) \in \Phi)$.
\begin{pro} \label{pro:compcong}
     Let $\ob{L}$ be a bialgebraic modular lattice, and let $a, b \in \ob{L}$ with
     $a \lcover b$. 
      Let 
       \( \Phi := \{ (x,y) \in \ob{L} \times \ob{L} \setsuchthat
                     \Int{a}{b} \text{ does not project into } \Int{x \meet y}{x \join y} \}.
       \)
     Then $\Phi$ is a complete congruence on the lattice $\ob{L}$.
\end{pro}
    Reflexivity and symmetry of $\Phi$ are obvious.  
    For transitivity, we assume that $(x, y) \in \Phi$, $(y, z) \in \Phi$ and
    $(x, z) \not\in \Phi$. 
    Since $(x, z) \not\in \Phi$, $\Int{a}{b}$ projects into $\Int{x \meet z}{x \join z}$, and hence
    into $\Int{x}{x \join z}$ or into $\Int{z}{x \join z}$.
    We will now distinguish two cases:
    \begin{enumerate}
        \item We assume that $\Int{a}{b}$ projects into 
                             $\Int{x}{x \join z}$: Then 
                             let $a_1, b_1 \in \ob{L}$ be such that
                             $x \le a_1 < b_1 \le x \join z$ and $\Int{a}{b} \projective \Int{a_1}{b_1}$.
                             We choose $\eta \in M(\ob{L})$ such that
              $\Int{a_1}{b_1} \nearrow \Int{\eta}{\eta^+}$.
            \begin{enumerate}
                 \item \emph{Case} $\eta \not\ge y$:  Then $\Int{\eta}{\eta^+} \searrow
                               \Int{\eta \meet y}{\eta^+ \meet y}$, and therefore
                               $\Int{a}{b}$ projects into $\Int{x \meet y}{y}$, implying
                               $(x, y) \not\in \Phi$, a contradiction.
                 \item \emph{Case} $\eta \ge y$: Since $\eta \ge x$ and $\eta \not\ge x \join z$,
                               we have $\eta \not\ge z$, and therefore
                              $\Int{\eta}{\eta^+} \searrow \Int{\eta \meet z}{\eta^+ \meet z}$. Hence
                               $\Int{a}{b}$ projects into $\Int{y \meet z}{z}$, implying
                               $(y, z) \not\in \Phi$, a contradiction.
            \end{enumerate}
       \item We assume that $\Int{a}{b}$ projects into 
                             $\Int{z}{x \join z}$: Then swapping the roles of $x$ and $z$ in
             the previous case, we obtain that $\Int{a}{b}$ projects into $\Int{z \meet y}{y}$ or
             into $\Int{y \meet x}{x}$, again contradicting the assumptions.
   \end{enumerate}
   Next, we will prove that if $(x_i)_{i \in I}$ and $(y_i)_{i \in I}$ are families from
   $\ob{L}$ such that for all $i \in I$, $(x_i, y_i) \in \Phi$, we have
   $(\bigjoin_{i \in I} x_i, \bigjoin_{i \in I} y_i) \in \Phi$.
   Let $X := \bigjoin_{i \in I} x_i$ and $Y :=  \bigjoin_{i \in I} y_i$.
   Seeking a contradiction, we assume $(X, Y) \not\in \Phi$. 
   Then $\Int{a}{b}$ projects into $\Int{X \meet Y}{X \join Y}$, and hence
   into $\Int{X}{X \join Y}$ or into $\Int{Y}{X \join Y}$. 
   In the case that $\Int{a}{b}$ projects into $\Int{X}{X \join Y}$, we choose
   $a_1, b_1 \in \ob{L}$ with $X \le a_1 < b_1 \le X \join Y$ and 
   $\Int{a}{b} \projective \Int{a_1}{b_1}$. We pick $\eta \in M(\ob{L})$
   with $\Int{a_1}{b_1} \nearrow \Int{\eta}{\eta^+}$. Since $\eta \not\ge b_1$,
   we have $\eta \not\ge X \join \bigjoin_{i \in I} y_i$. Since $\eta \ge a_1 \ge X$,
   there is $j \in I$ such that $\eta \not\ge y_j$. Then
   $\Int{\eta}{\eta^+} \searrow \Int{\eta \meet y_j}{\eta^+ \meet y_j}$,
   and therefore $\Int{a}{b}$ projects into $\Int{x_j \meet y_j}{y_j}$,
   implying $(x_j, y_j) \not\in \Phi$, a contradiction.
   In the case that $\Int{a}{b}$ projects into $\Int{Y}{X \join Y}$,
   we swap the roles of $X$ and $Y$ and obtain that $\Int{a}{b}$ projects
    into some $\Int{y_j \meet x_j}{x_j}$.
  Hence $\Phi$ is preserved under arbitrary joins.
    
   Now let $\ob{K}$ be the dual of $\ob{L}$, and let 
   $\Psi :=
              \{ (x,y) \in \ob{K} \times \ob{K} \setsuchthat
                     \Int{b}{a} \text{ does not project into } \Int{x \meet_{\ob{K}} y}{x \join_{\ob{K}} y} \text{ in }\ob{K} \}$.
   Since $\ob{L}$ is bialgebraic, so is $\ob{K}$, and hence the previous arguments
   imply that $\Psi$ is invariant under arbitrary joins, computed in $\ob{K}$.
   Hence $\Psi$ is invariant under arbitrary meets, computed in $\ob{L}$, and since
   $\Psi = \Phi$, we obtain that $\Phi$ is preserved under arbitary meets.
   
   Hence $\Phi$ is indeed a complete congruence of the lattice $\ob{L}$. \qed

We will also need some additional information on lonesome meet irreducible elements.
\begin{pro} \label{pro:lonesomeM}
    Let $\ob{L}$ be an algebraic modular lattice,
    and let $\eta \in M(\ob{L})$.
             If $\eta$ is not a lonesome meet irreducible
             element, then there exists $\psi \in M(\ob{L})$ with
             $\Int{\eta}{\eta^+} \projective \Int{\psi}{\psi^+}$,
             $\eta \not\le \psi$, and $\psi \not\le \eta$.
\end{pro}
\emph{Proof:}
     We let 
     $\varphi \in M(\ob{L})$ with $\eta \neq \varphi$ such that
    $\Int{\eta}{\eta^+} \projective
    \Int{\varphi}{\varphi^+}$.
    Since $\Int{\eta}{\eta^+} \projective \Int{\varphi}{\varphi^+}$,
     there is a natural number $n$, and there are $\rho_1, \ldots, \rho_{2n-1},
     \sigma_1, \ldots, \sigma_{2n-1} \in \ob{L}$ such that
       \begin{multline*}
        \Int{\eta}{\eta^+} \searrow
        \Int{\rho_1}{\sigma_1} \nearrow
        \Int{\rho_2}{\sigma_2} \searrow \\
         \cdots  \nearrow
        \Int{\rho_{2n - 2}}{\sigma_{2n-2}}
        \searrow
            \Int{\rho_{2n-1}}{\sigma_{2n-1}} 
              \nearrow
             \Int{\varphi}{\varphi^+}.
       \end{multline*}
       Now for each  $k \in \{1,2,\ldots,n-1\}$, 
       we pick
      an element $\eta_{2k} \in M (\ob{L})$ with
      $\eta_{2k} \geq \rho_{2k}$, $\eta_{2k} \not\geq \sigma_{2k}$.
      Then by Proposition~\ref{pro:projprops}\eqref{it:mi2},
      $\Int{\rho_{2k}}{\sigma_{2k}} \nearrow
       \Int{\eta_{2k}}{\eta_{2k}^+}$.
      Since $\nearrow$ is transitive, we obtain 
      \begin{multline}  \label{eq:chain2}
        \Int{\eta}{\eta^+} \searrow
        \Int{\rho_1}{\sigma_1} \nearrow
        \Int{\eta_2}{{\eta_2^+}} \searrow \\
       \cdots \nearrow
        \Int{\eta_{2n - 2}}{{\eta_{2n-2}^+}}
              \searrow
        \Int{\rho_{2n-1}}{\sigma_{2n-1}} 
              \nearrow
        \Int{\varphi}{\varphi^+}.
       \end{multline}
              Hence there exists $i \in \{1,\ldots, n\}$ and  $\psi$ with
         $\psi \neq \eta$ and 
         $\Int{\eta}{\eta^+} \searrow \Int{\rho_{2i-1}}{\sigma_{2i-1}} \nearrow \Int{\psi}{\psi^+}$.
         If $\eta \le \psi$, then $\eta^+ \le \psi$, and therefore
         $\rho_{2i-1} =
          \psi \meet \sigma_{2i-1} \ge \eta^+ \meet \sigma_{2i-1} = \sigma_{2i-1}$, 
         which implies $\psi^+  = \psi \join \sigma_{2i-1} \le \psi \join \rho_{2i-1} 
         = \psi$, a contradiction. 
         Hence $\eta \not\le \psi$. Similarly, we obtain $\psi \not\le \eta$. \qed

\begin{pro} \label{pro:lonesomeJ} 
   Let $\ob{L}$ be a bialgebraic modular lattice,
   and let $\beta \in J(\ob{L})$ and $\eta \in M(\ob{L})$ such that         
 $\Int{\beta^-}{\beta} \projective \Int{\eta}{\eta^+}$.
   Then    
     $\beta$ is a lonesome join irreducible element of $\ob{L}$
    if and only if $\eta$ is a lonesome meet irreducible element
    of $\ob{L}$.
\end{pro}
      For proving the ``only if''-direction, we assume
      that $\eta$ is a not a lonesome meet irreducible element of 
      $\ob{L}$. Let $\eta_1 \in M (\ob{L})$ be such that
      $\eta_1 \ge \beta^-$, $\eta_1 \not\ge \beta$.
       Then $\Int{\beta^-}{\beta} \nearrow \Int{\eta_1}{\eta_1^+}$.
     Since $\Int{\eta_1}{\eta_1^+} \projective \Int{\eta}{\eta^+}$, $\eta_1$ is not a lonesome meet irreducible
      element of $\ob{L}$, either.
      Let $\psi$ be a meet irreducible element with
      $\Int{\eta_1}{\eta_1^+} \projective \Int{\psi}{\psi^+}$ and $\eta_1 \not\le \psi$
       as
       produced in
      Proposition~\ref{pro:lonesomeM}. Using that the dual of $\ob{L}$ is algebraic
      and the dual of
      Proposition~\ref{pro:projprops}\eqref{it:mi3},
      we can choose  $\epsi \in J(\ob{L})$ such that $\epsi \le \eta_1$ and
      $\Int{\epsi^-}{\epsi} \nearrow \Int{\psi}{\psi^+}$.
      Since $\epsi \le \eta_1$ and $\beta \not\le \eta_1$,
      we have $\beta \neq \epsi$. Therefore, $\beta$ is not lonesome.
     
          The ``if''-direction now follows by applying the direction
         that has already been proved to the dual of $\ob{L}$.
          \qed 
  
    We notice that for finite lattices Proposition~\ref{pro:lonesomeJ} also follows from Corollary~6.2.1 
    of \cite{Av:DSOP}. To see this, we let $Q$ be
    the set of all prime intervals in $\ob{L}$ that are projective to $\Int{\beta^-}{\beta}$,
    and use Corollary~6.2.1 to establish that $\Int{\eta}{\eta^+}$ is the
    only element $\Int{x}{y}$ of $Q$ where $x$ is meet irreducible.

For a prime interval $\Int{\alpha}{\beta}$ with $\mygamma (\alpha, \beta) \neq 1_A$ 
in a bialgebraic modular lattice,
we will find a splitting of the congruence lattice. 
For any complete lattice $\ob{L}$ and $\alpha, \beta \in \ob{L}$, we define 
\[ 
    \begin{array}{rcl}
         \mydelta (\alpha, \beta) & := & \bigjoin \{ \rho \in J (\ob{L}) \setsuchthat
                                     \Int{\rho^-}{\rho} \projective \Int{\alpha}{\beta} \}.
    \end{array}
\]
\begin{pro} \label{pro:splits}
    Let $\ob{L}$ be an algebraic modular lattice, 
    let $\alpha, \beta \in \ob{L}$ with $\alpha \lcover \beta$.
     Then for all $\varphi \in \ob{L}$, we have
    $\varphi \le \mygamma (\alpha, \beta)$ or $\varphi \ge \mydelta (\alpha, \beta)$.
\end{pro}
\emph{Proof:} 
   Assume $\varphi \not\ge \mydelta (\alpha, \beta)$. By the definition of $\mydelta (\alpha, \beta)$,
     there is $\rho \in J(\ob{L})$ such that
   $\Int{\rho^-}{\rho} \projective \Int{\alpha}{\beta}$ and $\rho \not\le \varphi$.
     Using~Proposition~\ref{pro:projprops}\eqref{it:mi3}, we find a $\psi \in M(\ob{L})$ such that
    $\varphi \le \psi$ and $\Int{\rho^-}{\rho} \nearrow \Int{\psi}{\psi^+}$.
      By the definition of $\mygamma (\alpha, \beta)$, we have $\psi \le \mygamma (\alpha, \beta)$, and therefore
     $\varphi \le \psi \le \mygamma (\alpha, \beta)$. \qed

In a bialgebraic modular lattice, we can describe lonesome meet irreducible elements.
\begin{pro} \label{pro:charlonesome}
   Let $\ob{L}$ be a bialgebraic modular lattice, let $\alpha, \beta \in \ob{L}$ with $\alpha \lcover \beta$,
   and let $\eta \in M(\ob{L})$ with $\Int{\alpha}{\beta} \projective \Int{\eta}{\eta^+}$.
   Then the following are equivalent.
   \begin{enumerate}
   \item \label{it:e1} $\eta$ is lonesome.
   \item \label{it:e2} $\eta$ is completely meet prime.
   \item \label{it:e3} There is a complete lattice homomorphism $h : \ob{L} \to \ob{B}_2$ with
         $h(\alpha) = h (\eta) = 0$ and $h(\beta) = h (\eta^+) = 1$.
   \end{enumerate}
\end{pro}
\emph{Proof:}
   \eqref{it:e1}$\Rightarrow$\eqref{it:e2}:
     We assume that $\eta$ is lonesome. Then $\mygamma (\eta, \eta^+) = \eta$. 
     We choose $\rho \in J(\ob{L})$ with $\rho \le \eta^+$, $\rho \not\le \eta$.
     Then $\Int{\rho^-}{\rho} \nearrow \Int{\eta}{\eta^+}$, and therefore
     $\mydelta (\eta, \eta^+) \ge \rho$. Since $\eta \not\ge \rho$, this implies
      $\eta \not\ge \mydelta (\eta, \eta^+)$.
     Let $X \subseteq \ob{L}$ be such that
     $\bigmeet X \le \eta$. Seeking a contradiction, we assume that
     for all $x \in X$, we have $x \not\le \eta$. Then by Proposition~\ref{pro:splits},
     we obtain $x \ge \mydelta(\eta, \eta^+)$ for all $x \in X$, and therefore
     $\eta \ge \bigmeet X \ge \mydelta(\eta, \eta^+)$, a contradiction.
   \eqref{it:e2}$\Rightarrow$\eqref{it:e3}: 
      For $x \in \ob{L}$, we define
       $h(x) = 0$ if $x \le \eta$ and $h (x) = 1$ if $x \not\le \eta$. 
      Let $\theta \in \ob{L}$ be defined by $\theta = \bigmeet \{y \in \ob{L} \setsuchthat
       y \not\le \eta \}$. Since $\eta$ is completely meet prime, $\theta \not\le \eta$. 
       Thus for all $x \in \ob{L}$, we have $x \not\le \eta$ if and only if
       $x \ge \theta$.
       Hence $h (x) = 1$ if and only if $x \ge \theta$.
       Thus if $h (\bigjoin_{i \in I} x_i) = 0$, then
       $\bigjoin_{i \in I} x_i \le \eta$, and therefore for each $i \in I$,
      $h(x_i) = 0$, implying $\bigjoin_{i \in I} h( x_i ) = 0$, and if
       $h (\bigjoin_{i \in I} x_i) = 1$, then $\bigjoin_{i \in I} x_i \not\le \eta$,
       hence there is $j \in I$ with $x_j \not\le \eta$, and thus
       $\bigjoin_{i \in I} h (x_i) \ge h (x_j) = 1$. 
       Furthermore, if $\bigmeet_{i \in I} h (x_i) = 1$, then
       for all $i \in I$, we have $h(x_i) = 1$ and thus $x_i \ge \theta$. 
        Hence
       $\bigmeet_{i \in I} x_i \ge \theta$, and therefore
        $h (\bigmeet_{i \in I} x_i) = 1$. This is the essential
       step in proving that 
        $h$ is also
       a complete meet homomorphism.
       Now $h (\eta) = 0$ and $h (\eta^+) = 1$. Since $\alpha \le \beta$,
       we have $h(\alpha)  \le h(\beta)$. If $h (\alpha) = h (\beta)$,
       then $(\alpha, \beta)$ lies in the congruence $\mathrm{ker} (h)$.
       Since $\Int{\eta}{\eta^+} \projective \Int{\alpha}{\beta}$, 
       $(\eta, \eta^+)$ lies in the congruence generated by $(\alpha, \beta)$,
       and thus $(\eta, \eta^+) \in \mathrm{ker} (h)$. This implies $h(\eta) = h (\eta^+)$,
      a contradiction. Thus $h (\alpha) < h (\beta)$, which implies $h (\alpha) = 0$ and $h (\beta) = 1$.
    \eqref{it:e3}$\Rightarrow$\eqref{it:e1}: 
        If $\eta$ is not lonesome, then by Proposition~\ref{pro:lonesomeM}, there
        is $\psi \in M (\ob{L})$ such that $\Int{\eta}{\eta^+} \projective \Int{\psi}{\psi^+}$
        and $\psi \not \le \eta$. Since $(\eta, \eta^+) \not\in \mathrm{ker} (h)$, we 
        have $(\psi, \psi^+) \not\in \mathrm{ker} (h)$, and therefore $h(\psi) = 0$.
        Thus $0 = h (\psi) \join h (\eta) = h (\psi \join \eta) \ge h (\eta^+) = 1$, a contradiction. \qed

\section{Commutator lattices} \label{sec:cl}

\subsection{Special elements in commutator lattices}

  \begin{lem} \label{lem:3lonesome}
    Let $(\ob{L}, \join, \meet, [.,.])$ be a commutator lattice
    and let $(.:.)$ its associated residuation.
    \begin{enumerate}
         \item  \label{it:31} If $\eta \in M (\ob{L})$ is such that $\eta = (\eta : \eta^+)$, then
                $\eta$ is lonesome. 
         \item  \label{it:32} If $\alpha \in J(\ob{L})$ is such that
   $[\alpha, \alpha] = \alpha$, then $\alpha$ is lonesome.
     \end{enumerate}
\end{lem}
\emph{Proof:}
     \eqref{it:31} We assume that $\psi \in M(\ob{L})$ is such that
                   $\Int{\eta}{\eta^+} \projective \Int{\psi}{\psi^+}$.
                  Since $\eta^+ \not\le (\eta : \eta^+)$, we have
                  $[\eta^+, \eta^+] \not\le \eta$, and therefore by
                  Lemma~\ref{lem:centproj}, $[\psi^+, \psi^+] \not\le \psi$.
                  Therefore $(\psi : \psi^+) = \psi$. Using Lemma~\ref{lem:centproj}
                  again, we obtain $\eta = (\eta : \eta^+) = (\psi : \psi^+) = \psi$. 

      \eqref{it:32} We assume that $\beta \in J(\ob{L})$ is such that
                   $\Int{\alpha^-}{\alpha} \projective \Int{\beta^-}{\beta}$.
                    Since $[\alpha, \alpha] \not\le \alpha^-$, Lemma~\ref{lem:centproj}
                    yields $[\beta, \beta] \not\le \beta^-$. Therefore $[\beta, \beta] = \beta$.
                    Since $(\beta^- : \beta) \not\ge \beta$, Lemma~\ref{lem:centproj}
                    yields $(\alpha^- : \alpha) \not\ge \beta$, which implies
                    $[\beta, \alpha] \not\le \alpha^-$. Thus $[\beta, \alpha] = \alpha$,
                    and therefore $\alpha \le \beta$. 
                    Exchanging $\alpha$ and $\beta$, we obtain $\beta \le \alpha$, and
                    therefore $\alpha = \beta$. \qed

\subsection{Constructions of commutator operations}

In this section, we will provide three constructions of commutator multiplications on 
a given lattice $\ob{L}$.
For a complete lattice $\ob{K}$, a complete sublattice $\ob{L}$ of $\ob{K}$, and
an element $x \in \ob{K}$, we define its $\ob{L}$-closure $c_{\ob{L}} (x) :=
\bigmeet\{ y \in \ob{L} \setsuchthat y \ge x \}$.
This operation is $c_{\ob{L}}$ is a monotonic operation from $\ob{K}$ to $\ob{L}$, and
$c_{\ob{L}} (x) \ge x$ for  all $x \in \ob{K}$.
 \begin{pro} \label{pro:ck} 
    Let $\ob{K}$ be a complete lattice, and let $\ob{L}$ be a complete sublattice
    of $\ob{L}$. Then for all families $(x_i)_{i \in I}$ from $\ob{K}$,
    we have
    $c_{\ob{L}} ( \bigjoin_{i \in I} x_i ) = \bigjoin_{i \in I} c_{\ob{L}} (x_i)$.
\end{pro}
\emph{Proof:} 
   For $\ge$, we let $j \in I$. Then $c_{\ob{L}} (\bigjoin_{i \in I} x_i) \ge
   \bigjoin_{i \in I} x_i \ge x_j$, and thus 
    $c_{\ob{L}} (\bigjoin_{i \in I} x_i) \ge c_{\ob{L}} (x_j)$. Hence
    $c_{\ob{L}} (\bigjoin_{i \in I} x_i) \ge \bigjoin_{j \in I} c_{\ob{L}} (x_j)$.
    For $\le$, we observe that $\bigjoin_{i \in I} c_{\ob{L}} (x_i) \ge
                                \bigjoin_{i \in I} x_i$, and thus
    $\bigjoin_{i \in I} c_{\ob{L}} (x_i) \ge c_{\ob{L}} (\bigjoin_{i \in I} x_i)$. \qed

\begin{lem} \label{lem:c1}
    Let $(\ob{K}, \join, \meet, [.,.]_{\ob{K}})$ be a commutator lattice,
    and let $\ob{L}$ be a complete sublattice of $\ob{K}$. For $x, y \in \ob{L}$,
    we define
    $[x, y]_{\ob{L}} := c_{\ob{L}} ([x, y]_{\ob{K}})$. 
    Then
    $(\ob{L}, \join, \meet, [.,.]_{\ob{L}})$ is a commutator lattice, and
    we have $[x, y]_{\ob{L}} \ge [x, y]_{\ob{K}}$ for all $x, y \in \ob{L}$.
\end{lem}
\emph{Proof:}
     In order to show that $[.,.]_{\ob{L}}$ is a commutator multiplication on $\ob{L}$,
     we observe that for all $x, y \in \ob{L}$, we have
     $[x,y]_{\ob{L}} = c_{\ob{L}} ([x,y]_{\ob{K}}) = c_{\ob{L}} ([y,x]_{\ob{K}}) = [y,x]_{\ob{L}}$.
     Since $x \meet y  \in \ob{L}$ and $[x,y]_{\ob{K}} \le x \meet y$, we also have
     $c_{\ob{L}} ([x,y]_{\ob{K}}) \le x \meet y$, and therefore $[x,y]_{\ob{L}} \le x \meet y$.
     Now let $(x_i)_{i \in I}$ be a family from $\ob{L}$. Then
     Proposition~\ref{pro:ck} yields
     $[\bigjoin_{i \in I} x_i, y]_{\ob{L}} =
      c_{\ob{L}} ([\bigjoin_{i \in I} x_i, y]_{\ob{K}}) =
      c_{\ob{L}} (\bigjoin_{i \in I} [x_i, y]_{\ob{K}}) =
      \bigjoin_{i \in I} c_{\ob{L}} ([x_i, y]_{\ob{K}}) =
      \bigjoin_{i \in I} [x_i, y]_{\ob{L}}$.
     Hence $(\ob{L}, \join, \meet, [.,.]_{\ob{L}})$ is a commutator lattice. 
     Finally $[x,y]_{\ob{L}} = c_{\ob{L}} ([x,y]_{\ob{K}}) \ge [x,y]_{\ob{K}}$. \qed

\begin{lem} \label{lem:c2} 
    Let $\ob{L}$ be a complete lattice, and let $(\ob{K}, \join, \meet, [.,.]_{\ob{K}})$
    be a commutator lattice. We assume that $h$ is a complete lattice homomorphism
    from $\ob{L}$ to $\ob{K}$. For $x, y \in \ob{L}$, we define
    $[x, y]_{\ob{L}} := \bigmeet \{z \in \ob{L} \setsuchthat h (z) \ge [h (x), h (y)]_{\ob{K}} \}$.
    Then
    $(\ob{L}, \join, \meet, [.,.]_{\ob{L}})$ is a commutator lattice, and
    we have $h ([x, y]_{\ob{L}}) \ge [h(x),h(y)]_{\ob{K}}$ for all $x, y \in \ob{L}$.
\end{lem}
   \emph{Proof:}
     We fix $x, y \in \ob{L}$.
     For commutativity,  we observe that
     $[x, y]_{\ob{L}} = \bigmeet \{z \in \ob{L} \setsuchthat h (z) \ge [h (x), h (y)]_{\ob{K}} \} =
        \bigmeet \{z \in \ob{L} \setsuchthat h (z) \ge [h (y), h (x)]_{\ob{K}} \} = [y, x]_{\ob{L}}$.
     Since $h (x \meet y)  = h(x) \meet h(y) \ge [h (x), h(y)]_{\ob{K}}$, we have
     $x \meet y \ge [x,y]_{\ob{L}}$. Furthermore,
     $h ([x, y]_{\ob{L}}) =
             h (\bigmeet \{z \in \ob{L} \setsuchthat h (z) \ge [h (x), h (y)]_{\ob{K}} \}) =
             \bigmeet \{ h(z) \setsuchthat z \in \ob{L}, h (z) \ge [h (x), h (y)]_{\ob{K}} \}
             \ge [h(x), h(y)]_{\ob{K}}$. What remains to show is join distributivity.
      Let $(x_i)_{i \in I}$ be a family from $\ob{L}$, and let
      $j \in I$. Then $h ( [ \bigjoin_{i \in I} x_i, y ]_{\ob{L}} )  \ge
                       [ h(\bigjoin_{i \in I} x_i)  , h(y)]_{\ob{K}} =
                       [ \bigjoin_{i \in I} h(x_i)  , h(y)]_{\ob{K}} =
                       \bigjoin_{i \in I} [h(x_i), h(y)]_{\ob{K}} \ge 
                       [h(x_j), h(y)]_{\ob{K}}$.
        Now by the definition of $[x_j, y]_{\ob{L}}$, this inequality implies
        $[ \bigjoin_{i \in I} x_i, y ]_{\ob{L}} \ge
                                       [ x_j, y ]_{\ob{L}}$. Therefore,
       $ [ \bigjoin_{i \in I} x_i, y ]_{\ob{L}} \ge \bigjoin_{i \in I} [ x_i, y]_{\ob{L}}$.
       For the other inequality, we observe that $h (\bigjoin_{i \in I} [x_i, y]_{\ob{L}}) =
       \bigjoin_{i \in I} h ([x_i, y]_{\ob{L}}) \ge \bigjoin_{i \in I} [h(x_i), h(y)]_{\ob{K}}
       = [\bigjoin_{i \in I} h (x_i), h(y)]_{\ob{K}} 
       = [h (\bigjoin_{i \in I} x_i), h(y)]_{\ob{K}}$. Using the definition
       of $[\bigjoin_{i \in I} x_i, y]_{\ob{L}}$, we obtain
        $\bigjoin_{i \in I} [x_i, y]_{\ob{L}} \ge 
         [\bigjoin_{i \in I} x_i, y]_{\ob{L}}$. \qed

\begin{pro} \label{lem:sisjoinhomo}
  Let $\ob{L}$ be a complete lattice, let $\Theta$ be a complete congruence
  on $\ob{L}$, and let $s :\ob{L} \to \ob{L}$ be the mapping defined by
  $s(x) := \bigmeet \{ z \in \ob{L} \setsuchthat (z, x) \in \Theta \}$.
  Then for all $x \in \ob{L}$, we have $(s(x), x) \in \Theta$, and
  for all families $(x_i)_{i \in I}$ from $\ob{L}$, we have
  $s (\bigjoin_{i \in I} x_i) = \bigjoin_{i \in I} s (x_i)$.
\end{pro}

\emph{Proof:}
  The fact that $\Theta$ is a complete congruence implies
  that for every $x \in \ob{L}$, we have
  $\bigmeet \{ z \in \ob{L} \setsuchthat (z, x) \in \Theta \}
   \equiv_{\Theta} x$, and hence 
  $(s (x), x) \in \Theta$.
  We first prove $s (\bigjoin_{i \in I} x_i) \ge \bigjoin_{i \in I} s (x_i)$. 
  To this end, let $j \in I$. We have 
  $s (\bigjoin_{i \in I} x_i) \equiv_{\Theta} \bigjoin_{i \in I} x_i$, and therefore
  $s (\bigjoin_{i \in I} x_i) \meet x_j \equiv_{\Theta} x_j$.
  From the definition of $s(x_j)$, we obtain
  $s(x_j) \le  s (\bigjoin_{i \in I} x_i) \meet x_j$, which implies
  $s (x_j) \le s (\bigjoin_{i \in I} x_i)$. Thus $\bigjoin_{i \in I} s (x_i) \le 
  s (\bigjoin_{i \in I} x_i)$. 
   For proving $s (\bigjoin_{i \in I} x_i) \le \bigjoin_{i \in I} s (x_i)$,
  we notice that $\bigjoin_{i \in I} s(x_i) \equiv_{\Theta} \bigjoin_{i \in I} x_i$.
  Hence from the definition of $s (\bigjoin_{i \in I} x_i)$, we obtain
  $s (\bigjoin_{i \in I} x_i) \le \bigjoin_{i \in I} s(x_i)$. \qed

\begin{lem} \label{lem:c3}
  Let $\ob{L}$ be a complete lattice that splits with splitting pair
  $(\delta, \epsi)$. Let $\Theta$ be a complete congruence of $\ob{L}$
  with $(\epsi, 1) \in \Theta$, and let $s$ be the complete join homomorphism
  associated with $\Theta$ that was defined in Lemma~\ref{lem:sisjoinhomo}.
  For $x, y \in \ob{L}$, we define
  $[x,y] := 0$ if $x \le \delta$ and $y \le \delta$, and
  $[x,y] := s (x \meet y) = \bigmeet \{ z \in \ob{L} : (z, x \meet y) \in \Theta \}$ otherwise.
   Then we have:
   \begin{enumerate}
     \item \label{it:c1} If $x \not\le \delta$, then $[x, y] = [y, x] = s (y)$.
     \item \label{it:c2} $[.,.]$ is a commutator multiplication on $\algop{\ob{L}}{\join, \meet}$.
   \end{enumerate}
\end{lem}
   
\emph{Proof:}
   For item~\eqref{it:c1}, we fix $x, y \in \ob{L}$ with $x \not\le \delta$. Then
   $x \ge \epsi$, and therefore $(x, 1) \in \Theta$.
   Then $[x, y] = s (x \meet y) = s (1 \meet y) = s(y) = s(y \meet 1) = s (y \meet x)
         = [y, x]$.
   For item~\eqref{it:c2}, we observe that
   commutativity and $[x,y] \le x \meet y$ for all $x,y \in \ob{L}$ follow
   immediately from the definition. 
   What remains to be proved is the join distributivity 
   $[\bigjoin_{i \in I} x_i, y] = \bigjoin_{i \in I} [x_i, y]$.
   In the case $y \not\le \delta$, item~\eqref{it:c1} yields
   $[\bigjoin_{i \in I} x_i, y] = 
    s (\bigjoin_{i \in I} x_i)$. By Lemma~\ref{lem:sisjoinhomo}, this last expression
   is equal to     
    $\bigjoin_{i \in I} s (x_i)$. Applying item~\eqref{it:c1} again,
   thie last expression is equal to $\bigjoin_{i \in I} [x_i, y]$.
    Next, we consider the case that $y \le \delta$ and
    $\bigjoin_{i \in I} x_i \le \delta$. 
    Then $[\bigjoin_{i \in I} x_i, y] = 0 = \bigjoin_{i \in I} 0 =
          \bigjoin_{i \in I} [x_i, y]$.
    The last case is that $y \le \delta$ and there exists $j \in I$ with
     $x_j \not\le \delta$. Then by item~\eqref{it:c1}, $[\bigjoin_{i \in I} x_i, y] = s (y)$.
    Now we compute $\bigjoin_{i \in I} [x_i, y] =
                    [x_j, y] \join \bigjoin_{i \in I \setminus \{j\}} [x_i, y]$.
    Again by~\eqref{it:c1}, the first joinand $[x_j, y]$ is equal to $s (y)$.
    For an arbitrary $i \in I \setminus \{j\}$, $[x_i, y] = 0$ if $x_i \le \delta$,
    and $s(y)$ if $x_i \not\le \delta$.
    Thus we have  $\bigjoin_{i \in I \setminus \{j\}} [x_i, y] \le s (y)$.
    Hence $\bigjoin_{i \in I} [x_i, y] = s(y)$, and so the equation expressing join distributivity 
    also holds in this last case. \qed 

 This construction had one origin in the analysis of \cite[Proposition~16]{IS:PRA}.

\subsection{Types of commutator lattices}

 \begin{de}
   Let $\ab{L} = \algop{\ob{L}}{\join, \meet, [.,.]}$ be a commutator lattice. Let $(\gamma_n)_{n \in \N}$ and $(\lambda_n)_{n \in \N}$ be the sequences in $\ob{L}$
   defined by $\gamma_1 = \lambda_1 = 1$ and 
   $\gamma_{n+1} = [\gamma_n, \gamma_n]$ and $\lambda_{n+1} = [1, \lambda_n]$ for $n \in \N$.
   Then $\ab{L}$ is of \emph{solvable type} if there is an $n \in \N$ with
   $\gamma_n = 0$, of \emph{nilpotent type} if there is an $n \in \N$ with $\lambda_n = 0$, and
   of \emph{abelian type} if $\gamma_2 = 0$. 
 \end{de}

\begin{lem} \label{lem:solnil}
   Let $\ab{L} = (\ob{L}, \join, \meet, [.,.])$ be a commutator lattice of finite height, and let
   $(.:.)$ its associated residuation. Then we have:
   \begin{enumerate}
      \item \label{it:cs}  $\ab{L}$ is of solvable type if and only if there is no $\beta \in \ob{L}$ with  
            $\beta \neq 0$ and $[\beta, \beta] = \beta$.
      \item \label{it:cn1}  $\ab{L}$ is of nilpotent type if and only if there is no $\beta \in \ob{L}$ with  
            $\beta \neq 0$ and $[1_A, \beta] = \beta$.
      \item \label{it:cn2} Assume that $\ob{L}$ is modular.
            Then $\ab{L}$ is of nilpotent type if and only if for all $\alpha, \beta \in \ob{L}$
            with $\alpha \lcover \beta$, we have $(\alpha:\beta) = 1$.
   \end{enumerate}
\end{lem}
    \emph{Proof:}
    \eqref{it:cs}
   If $\ab{L}$ is not of solvable type, then there will be an $n \in \N$ such that
   in the derived series of $\ab{L}$ we have $\gamma_n = \gamma_{n+1}$ and $\gamma_n \neq 0$.
   Then we set $\beta := \gamma_n$. On the other hand, if $[\beta, \beta] = \beta$, we 
   prove by induction that $\gamma_n \ge \beta$ for all $n \in \N$.
   The induction step is $\gamma_{n+1} = [\gamma_n , \gamma_n] \ge [\beta, \beta] = \beta$.
   This proves item~\eqref{it:cs}. Item~\eqref{it:cn1} is proved similarly.
   For the ``if''-direction of~\eqref{it:cn2},
   we assume that
   $\ab{L}$ is not of nilpotent type. By~\eqref{it:cn1}, there is $\gamma \in \ob{L}$
   such that $\gamma \neq 0$ and $[1, \gamma] = \gamma$. Let $\delta \lcover \gamma$.
   Then $[1 , \gamma] \not\le \delta$, hence $(\delta : \gamma) \neq 1$.
   For the ``only if''-direction of~\eqref{it:cn2}, we assume that there are $\alpha \lcover \beta$
   in $\ob{L}$ such that $(\alpha : \beta) < 1$. Let $\gamma$ be minimal with
   $\gamma \le \beta$, $\gamma \not\le \alpha$. Then $\gamma$ is join irreducible and 
   $\Int{\gamma^-}{\gamma} \nearrow \Int{\alpha}{\beta}$. From Lemma~\ref{lem:centproj},
    we obtain $(\gamma^- : \gamma) =
    (\alpha:\beta) < 1$, and therefore $[1, \gamma] \not\le \gamma^-$. This implies
   $[1, \gamma] = \gamma$, and hence by item~\eqref{it:cn1}, $\ab{L}$ is not of nilpotent type. \qed

\subsection{The largest commutator operation on a given lattice}

Given a lattice $\ob{L}$  and $x,y \in \ob{L}$,
we would like to obtain an upper bound for $[x, y]$ for each
commutator multiplication definable on $\ob{L}$. From \cite{Cz:AOTC},
we know that such a bound is provided by the single largest commutator multiplication on
each lattice:
\begin{lem}[{\cite[Corollary~1.5]{Cz:AOTC}}]
  Let $\ob{L}$ be a complete lattice, and
  let $([.,.]_i)_{i \in I}$ be the family of all binary operations that
  turn $\ob{L}$ into a commutator lattice.
  For $x,y \in \ob{L}$, we define
  $\LComm{x}{y}_{\ob{L}} := \bigjoin_{i \in I} [x,y]_i$.
  Then $(\ob{L}, \join, \meet, \LComm{.}{.}_{\ob{L}})$ is a commutator lattice.
\end{lem}
Czelakowski writes $\bullet_{\Omega}$ for the operation $\LComm{.}{.}_{\ob{L}}$ 
and states that ``the characterization of the operation $\bullet_{\Omega}$ in modular
algebraic lattices is an open and challenging problem'' \cite[p.\ 114]{Cz:AOTC}.
We will not be able to construct this operation $\LComm{.}{.}_{\ob{L}}$ completely,
but we will obtain a description of the associated residuum $(\alpha:\beta)$
if $\alpha \lcover \beta$.

\begin{de}
  Let $\ob{L}$ be a complete lattice. Then $\ob{L}$ \emph{forces abelian type} if
  $\LComm{1}{1}_{\ob{L}} = 0$. $\ob{L}$ forces \emph{nilpotent type} if $(\ob{L}, \join, \meet, \LComm{.}{.}_{\ob{L}})$
  is of nilpotent type, and $\ob{L}$ forces \emph{solvable type} if
  $(\ob{L}, \join, \meet, \LComm{.}{.}_{\ob{L}})$ is of solvable type.
\end{de}

\begin{lem} \label{lem:klsub}
  Let $\ob{K}$ be a complete lattice, and let $\ob{L}$ be a complete sublattice of
  $\ob{K}$. Then for all $x,y \in \ob{L}$, we have
  $\LComm{x}{y}_{\ob{L}} \ge \LComm{x}{y}_{\ob{K}}$.
\end{lem}
\emph{Proof:}
    We use Lemma~\ref{lem:c1} to construct a multiplication $[.,.]_{\ob{L}}$ on $\ob{L}$
    by $[.,.]_{\ob{L}} := c_{\ob{L}} (\LComm{x}{y}_{\ob{K}})$ for $x, y \in \ob{L}$.
    Then for all $x, y \in \ob{L}$, we have 
    $\LComm{x}{y}_{\ob{L}} \ge [x, y]_{\ob{L}} = c_{\ob{L}} (\LComm{x}{y}_{\ob{K}}) \ge
     \LComm{x}{y}_{\ob{K}}$. \qed

\begin{lem} \label{lem:klhomo}
  Let $\ob{L}, \ob{K}$ be complete lattices, and let $h$ be a complete
  lattice homomorphism from $\ob{L}$ to $\ob{K}$.
  Then for all $x,y \in \ob{L}$, we have
  $h (\LComm{x}{y}_{\ob{L}}) \ge \LComm{h(x)}{h(y)}_{\ob{K}}$.
\end{lem}
    We use Lemma~\ref{lem:c2} to construct a multiplication $[.,.]_{\ob{L}}$ on $\ob{L}$
    by $[x, y]_{\ob{L}} := 
    \bigmeet \{ z \in \ob{L} \setsuchthat h(z) \ge \LComm{h(x)}{h(y)}_{\ob{K}} \}$
    for $x, y \in \ob{L}$.
    Now for all $x, y \in \ob{L}$, we have 
    $h (\LComm{x}{y}_{\ob{L}}) \ge h ([x, y]_{\ob{L}})$.
    By Lemma~\ref{lem:c2}, we have $h ([x,y ]_{\ob{L}}) \ge
         \LComm{x}{y}_{\ob{K}}$. \qed   

We call a complete sublattice $\ob{L}$ of $\ob{K}$ a \emph{complete $(0,1)$-sublattice}
of $\ob{K}$ 
if $0_{\ob{K}} \in \ob{L}$ and $1_{\ob{K}} \in \ob{L}$. In this case, $\ob{K}$ is
a \emph{$(0,1)$-extension} of $\ob{L}$.
\begin{thm} \label{thm:sub}
  Let $\ob{L}$ be a complete lattice, and let $\ob{K}$ be a complete $(0,1)$-extension or
  a complete $(0,1)$-homomorphic image of $\ob{L}$.
  If $\ob{L}$ forces abelian, nilpotent, or solvable type, then so does
  $\ob{K}$.
\end{thm}

 \emph{Proof:} We use a function $f$ to treat the lower central and the derived series
    at once.
    Let $f : \N \setminus \{1\} \to \N$ be a function with $f(n) < n$ for all $n \in \N$, and 
    let $(\kappa_n)_{n \in \N}$ be a sequence from $\ob{K}$ defined by
    $\kappa_1 = 1$ and $\kappa_n := \LComm{\kappa_{f(n)}}{\kappa_{n-1}}_{\ob{K}}$ for $n > 1$.
    Let $(\lambda_n)_{n \in \N}$ be the corresponding sequence from
    $\ob{L}$ defined by 
     $\lambda_1 = 1$ and $\lambda_n := \LComm{\lambda_{f(n)}}{\lambda_{n-1}}_{\ob{L}}$ for $n > 1$.

    If $\ob{L}$ is a $(0,1)$ sublattice of $\ob{K}$, we have
     $\kappa_n \le \lambda_n$ for all $n \in \N$. We prove this by induction:
     for $n > 1$, $\kappa_{n} = \LComm{\kappa_{f(n)}}{\kappa_{n-1}}_{\ob{K}}$. 
      By the induction hypothesis and monotonicity, we obtain
         $\LComm{\kappa_{f(n)}}{\kappa_{n-1}}_{\ob{K}}
         \le
         \LComm{\lambda_{f(n)}}{\lambda_{n-1}}_{\ob{K}}$. 
    By Lemma~\ref{lem:klsub}, we have $\LComm{\lambda_{f(n)}}{\lambda_{n-1}}_{\ob{K}} \le
    \LComm{\lambda_{f(n)}}{\lambda_{n-1}}_{\ob{L}} = \lambda_{n}$.
    Therefore, if for some $k \in \N$, $\lambda_k = 0$, then $\kappa_k = 0$.

    If $\ob{K}$ is a complete $(0,1)$-homomorphic image of $\ob{L}$, we have
    $\kappa_n \le h (\lambda_n)$ for all $n \in \N$. Again, we proceed by induction:
    as the induction basis, we observe that $\kappa_1 = 1_{\ob{K}} = h (1_{\ob{L}}) = h (\lambda_1)$.
    For the induction step, we let $n > 1$ and compute
      $\kappa_{n} = \LComm{\kappa_{f(n)}}{\kappa_{n-1}}_{\ob{K}}$. 
      By the induction hypothesis and monotonicity, we obtain
         $\LComm{\kappa_{f(n)}}{\kappa_{n-1}}_{\ob{K}}
         \le
         \LComm{h(\lambda_{f(n)})}{h(\lambda_{n-1})}_{\ob{K}}$. 
    By Lemma~\ref{lem:klhomo}, we have $\LComm{h(\lambda_{f(n)})}{h(\lambda_{n-1})}_{\ob{K}} \le
    h (\LComm{\lambda_{f(n)}}{\lambda_{n-1}}_{\ob{L}}) = h (\lambda_{n})$.
    Therefore, if for some $k \in \N$, $\lambda_k = 0$, then $\kappa_k = h (\lambda_k) = 0$.
     
     Now if $\ob{L}$ forces abelian type, then $\lambda_2 = 0$, and hence $\kappa_2 = 0$,
     and therefore $\ob{K}$ forces abelian type. If $\ob{L}$ forces nilpotent type,
     we choose $f(n) := 1$ for all $n \in \N$ and observe that there is $k \in \N$ with
     $\lambda_k = 0$, hence $\kappa_k = 0$, and thus $\ob{K}$ forces nilpotent type.
     For solvable type, the proof is analogous with $f(n) := n - 1$. \qed  
 
  \begin{thm} \label{thm:forceabel}
   Let $\ob{K}$ be a complete lattice. If $\ob{K}$ has a complete
   $(0,1)$-sublattice $\ob{L}$ that is algebraic, modular, simple, complemented, and has at least $3$ elements,
   then $\ob{K}$ forces abelian type.
 \end{thm}
 \emph{Proof:}
   By Theorem~\ref{thm:sub}, it is sufficient to prove that $\ob{L}$ forces abelian type.
    Let $T$ be the set of atoms of $\ob{L}$. We let $\lfloor .:. \rfloor$ denote the residuation operation
   associated with the largest commutator operation $\LComm{.}{.}_{\ob{L}}$ on $\ob{L}$.
   We show that for all $\alpha \in T$, $\LComm{\alpha}{\alpha}_{\ob{L}} = 0$.
   Let $\eta_1$ be a complement of $\alpha$ in $\ob{L}$. Then $\Int{0}{\alpha} \nearrow \Int{\eta_1}{1}$,
   and therefore $\eta_1$ is a coatom of $\ob{L}$.  Since $|\ob{L}| \ge 3$, $\eta_1 \neq 0$, and 
    therefore by \cite[Lemma~4.83]{MMT:ALVV}, there is an atom $\beta$ of $\ob{L}$ with
    $\beta \le \eta_1$. Let $\eta_2$ be the  complement of $\beta$ in $\ob{L}$. Then
    $\eta_2$ is a coatom of $\ob{L}$ and $\eta_1 \neq \eta_2$.
   By Lemma~\ref{lem:res}\eqref{it:l3}, we have $\lfloor \eta_1 : 1 \rfloor_{\ob{L}} \ge \eta_1$. Since $\ob{L}$ is
   simple and modular, Dilworth's congruence generation theorem
   \cite[Theorem~2.66]{MMT:ALVV} yields that the intervals $\Int{\eta_1}{1}$ and $\Int{\eta_2}{1}$ are projective inside $\ob{L}$.
   Hence by Lemma~\ref{lem:centproj},
   $\lfloor \eta_2 : 1 \rfloor_{\ob{L}}
   = \lfloor \eta_1 : 1 \rfloor_{\ob{L}} \ge \eta_1$. Since $\lfloor \eta_2 : 1 \rfloor_{\ob{L}} \ge \eta_2$, we obtain
   $\lfloor \eta_2 : 1 \rfloor_{\ob{L}} \ge \eta_1 \join \eta_2 =  1$. Thus $\LComm{1}{1}_{\ob{L}} \le \eta_2$, and therefore, since 
   by simplicity all
   prime intervals of $\ob{L}$ are projective, Lemma~\ref{lem:centproj} yields
   $\LComm{\alpha}{\alpha}_{\ob{L}} = 0$.
      In an algebraic complemented modular lattice, $1$ is the join of atoms \cite[Lemma~4.83]{MMT:ALVV}.
   Hence
   $\LComm{1}{1}_{\ob{L}} = \LComm{\bigjoin T}{\bigjoin T}_{\ob{L}} = \bigjoin \{ \LComm{\alpha}{\beta}_{\ob{L}} \setsuchthat \alpha, \beta \in T \}$.
   The joinands of the last expression with $\alpha = \beta$ are $0$ by the above argument;
   for the other joinands, we have $\LComm{\alpha}{\beta}_{\ob{L}} \le \alpha \meet \beta = 0$. 
   This completes the proof that $\ob{L}$ forces abelian type; now Theorem~\ref{thm:sub} implies
   that $\ob{K}$ forces abelian type. \qed

  \begin{thm} \label{thm:cab2}
    Let $\ob{L}$ be a 
    bialgebraic 
    modular lattice, and
    let $\lfloor x: y \rfloor :=
   \bigjoin \{ z \in \ob{L} \setsuchthat \LComm{z}{y}_{\ob{L}} \le x \}$  denote the residuation operation associated with
   $\LComm{.}{.}_{\ob{L}}$. Let $\alpha, \beta \in \ob{L}$ be such that
   $\alpha \lcover \beta$. Then 
   $\lfloor \alpha : \beta \rfloor = \mygamma (\alpha, \beta)$.
 \end{thm}
 
 \emph{Proof:} Lemma~\ref{lem:lb} implies $\lfloor \alpha : \beta \rfloor \ge
    \mygamma (\alpha, \beta)$. For proving $\le$, we let $\rho$ be a join irreducible
    element of $\ob{L}$ with $\rho \le \beta, \rho \not\le \alpha$. 
    Then $\Int{\rho^-}{\rho} \nearrow \Int{\alpha}{\beta}$, and therefore
    $\mygamma (\alpha, \beta) = \mygamma (\rho^-, \rho)$, and  
    by Lemma~\ref{lem:centproj}
    $\lfloor \alpha : \beta \rfloor = \lfloor \rho^- : \rho \rfloor$.
    In order to prove $\lfloor \rho^- : \rho \rfloor \le \mygamma (\rho^-, \rho)$,
    we fix $\psi \in \ob{L}$ such that 
    such that $\LComm{\psi}{\rho}_{\ob{L}} \le \rho^-$, and show that
    $\psi \le \mygamma(\rho^-, \rho)$.
    In the case $\mygamma (\rho^-, \rho) = 1$, this is obviously true, 
    so we
    assume $\mygamma(\rho^-, \rho) < 1$.
     We let $\delta := \mygamma(\rho^-, \rho)$ and 
     we define $\epsi \in \ob{L}$ by
      $\epsi :=  \mydelta (\rho^-, \rho)$, which is defined as $\bigjoin \{ \sigma \in J (\ob{L}) \setsuchthat
                                     \Int{\sigma^-}{\sigma} \projective \Int{\rho^-}{\rho} \}$.
     Then $(\delta, \epsi)$ is a splitting pair for the lattice $\ob{L}$.
     Let $\Theta$ be the complete congruence of $\ob{L}$ that is generated by
     $(\epsi, 1)$. Then we apply Lemma~\ref{lem:c3} to  $\Theta$ and the splitting pair $(\delta, \epsi)$
     and obtain a complete join homomorphism $s$ and a commutator multiplication $[.,.]$ on $\ob{L}$.
     Next, we show
     \begin{equation} \label{eq:1}
          (\rho^-, \rho) \not\in \Theta.
     \end{equation}
      For this purpose, we construct a complete congruence $\Phi$ of $\ob{L}$ such 
      that $(\epsi, 1) \in \Phi$ and $(\rho^-, \rho) \not\in \Phi$. 
      Let $a := \rho^-$ and $b := \rho$, let $\Phi$ be the complete congruence
      of $\ob{L}$ produced in Proposition~\ref{pro:compcong}. Clearly,
      $(\rho^-, \rho) \not\in \Phi$. 
      Now 
      suppose $(\epsi, 1) \not\in \Phi$. Then 
      we have  $\rho_1, \rho_2 \in \ob{L}$ such that
      $\epsi \le \rho_1 \lcover \rho_2$ and $\Int{\rho^-}{\rho} \projective
       \Int{\rho_1}{\rho_2}$. From the dual of Proposition~\ref{pro:projprops}\eqref{it:mi1},
       we obtain $\rho_3 \in J (\ob{L})$ with $\rho_3 \le \rho_2$, $\rho_3 \not\le \rho_1$.
       Then by the dual of Proposition~\ref{pro:projprops}\eqref{it:mi2}, 
      $\Int{\rho_3^-}{\rho_3} \nearrow \Int{\rho_1}{\rho_2} \projective \Int{\rho^-}{\rho}$, and
      therefore from the definition of $\epsi$ as
      $\mydelta (\rho^-, \rho)$, we obtain $\rho_3 \le \epsi$. Thus $\rho_3 \le \rho_1$, and therefore $\rho_2 = \rho_1 \join \rho_3
       = \rho_1$, a contradiction. This contradiction proves $(\epsi, 1) \in \Phi$.
       Hence $\Theta \subseteq \Phi$, which completes the proof of~\eqref{eq:1}.
      We will next prove $s(\rho) = \rho$. Suppose $s (\rho) < \rho$. Then 
      $s (\rho) \le \rho^-$. By Proposition~\ref{lem:sisjoinhomo},
      we have $(s(\rho), \rho) \in \Theta$, and therefore $(\rho^-, \rho) = (s (\rho) \join \rho^-, \rho \join \rho^-) \in \Theta$,
      contradicting~\eqref{eq:1}.
      Therefore $s(\rho) = \rho$.
      Since $\LComm{\psi}{\rho}_{\ob{L}} \le \rho^-$, we have $[\psi, \rho] \le \rho^-$.
      Now if $\psi \not\le \delta$, then by Lemma~\ref{lem:c3}\eqref{it:c1}, $\rho^- \ge [\psi, \rho] = s(\rho) = \rho$,
      a contradiction. Therefore $\psi \le \delta = \mygamma(\rho^-, \rho) = \mygamma (\alpha, \beta)$. \qed

 \begin{thm} \label{thm:forcenil}
   Let $\ob{L}$ be a modular lattice of finite height.
   Then $\ob{L}$ forces nilpotent type if and only if
   for all $\alpha, \beta \in \ob{L}$ with $\alpha \lcover \beta$, we
   have $\mygamma(\alpha, \beta) = 1$.
 \end{thm}
    We let $\lfloor x: y \rfloor :=
   \bigjoin \{ z \in \ob{L} \setsuchthat \LComm{z}{y}_{\ob{L}} \le x \}$  denote the residuation operation associated with
  $\LComm{.}{.}_{\ob{L}}$. In order to show that $(\ob{L}, \join, \meet, \LComm{.}{.}_{\ob{L}})$ is of nilpotent type,
   we use Lemma~\ref{lem:solnil}~\eqref{it:cn2}. 
   By this Lemma, $\ob{L}$ forces nilpotent type if and only if 
   for all $\alpha, \beta \in \ob{L}$ with  $\alpha \lcover \beta$,
   we have $\lfloor \alpha : \beta \rfloor = 1$, which    
   by Theorem~\ref{thm:cab2} is equivalent to $\mygamma (\alpha, \beta) = 1$.
   \qed

  Next, we want to characterize lattices forcing solvable type.
  \begin{lem}  \label{lem:lone}
 Let $\ob{L}$ be a bialgebraic 
   modular lattice, 
    let $\lfloor x: y \rfloor :=
   \bigjoin \{ z \in \ob{L} \setsuchthat \LComm{z}{y}_{\ob{L}} \le x \}$  denote the residuation operation associated with
   $\LComm{.}{.}_{\ob{L}}$, and let $\eta \in M( \ob{L} )$.
       Then $\lfloor \eta : \eta^+ \rfloor = \eta$ if and only if $\eta$ is lonesome.
   \end{lem}
   \emph{Proof:}
      The ``only if''-direction is a consequence of Lemma~\ref{lem:3lonesome}. 
      For the ``if''-direction, we assume that $\eta$ is lonesome. 
      Then Proposition~\ref{pro:charlonesome} yields a complete lattice homomorphism from $\ob{L}$ onto $\ob{B}_2$
      with $h (\eta) = 0$ and $h(\eta^+) = 1$. Now
      Lemma~\ref{lem:klhomo} implies $h (\LComm{\eta^+}{\eta^+}_{\ob{L}}) \ge \LComm{h (\eta^+)}{h (\eta^+)}_{\ob{B}_2} =
      \LComm{1}{1}_{\ob{B}_2}$,  which is equal to $1$
    because
    $[x, y] := x \meet y$ is a commutator multiplication on $\ob{B}_2$. Therefore $\LComm{\eta^+}{\eta^+}_{\ob{L}} \not\le \eta$,
      and then
      $\eta^+ \not\le \lfloor \eta : \eta^+ \rfloor$. Since $\eta \le \lfloor \eta : \eta^+ \rfloor$,
      we have $\lfloor \eta : \eta^+ \rfloor = \eta$. \qed

 \begin{thm} \label{thm:lattsolv}
   Let $\ob{L}$ be a modular lattice of finite height.
   Then $\ob{L}$ forces solvable type if and only if the two element lattice $\ob{B}_2$ is not
   a homomorphic image of $\ob{L}$.
 \end{thm}
 \emph{Proof:}
    For the ``only if''-direction, we assume that $\ob{L}$ forces solvable type and that
    $h : \ob{L} \to \ob{B}_2$ is an epimorphism.
    Then by Theorem~\ref{thm:sub}, $\ob{B}_2$ forces solvable type, which contradicts
    the fact that on $\ob{B}_2$, the operation $[x,y] := x \meet y$ is a commutator
    multiplication which is not of solvable type. 
    For the ``if''-direction, we assume that $\ob{L}$ does not force solvable type.
    Then by Lemma~\ref{lem:solnil}, there is a $\beta \in \ob{L}$ with $\beta > 0$ and
    $\LComm{\beta}{\beta}_{\ob{L}} = \beta$.
     Let $\alpha \lcover \beta$, and let $\rho$ be minimal with $\rho \le \beta$, $\rho \not\le \alpha$.
     Then $\rho$ is join irreducible and $\LComm{\rho}{\rho}_{\ob{L}} \not\le \rho^-$, and therefore
     $\LComm{\rho}{\rho}_{\ob{L}} = \rho$. Now by Lemma~\ref{lem:3lonesome}, $\rho$ is lonesome.
    Taking $\eta \in M(\ob{L})$ with $\eta \ge \rho^-$ and  $\eta \not\ge \rho$
    and using Propositions~\ref{pro:projprops} and \ref{pro:lonesomeJ} we obtain that $\eta$ is lonesome, 
    and now Proposition~\ref{pro:charlonesome} yields an epimorphism of $\ob{L}$ onto
    $\ob{B}_2$. \qed 
   
\section{Algebras}

\subsection{Lattice conditions}

The results on commutator multiplications of the previous sections 
immediately yield the following theorem.
\begin{thm} \label{thm:lattforce}
   Let $\ab{A}$ be an algebra in a congruence modular variety.
   \begin{enumerate}
       \item \label{it:tf1} If $\Con (\ab{A})$ has a complete $(0,1)$-sublattice            with at least $3$ elements that
             is algebraic, 
             simple, and complemented,
             then $\ab{A}$ is abelian.
      \item \label{it:tf3} If $\Con (\ab{A})$ has a finite $(0,1)$-sublattice $\ob{L}$ that
        does not split, then $\ab{A}$ is supernilpotent.
         \item \label{it:tf2} If $\Con (\ab{A})$ has a $(0,1)$-sublattice $\ob{L}$ of finite height 
             such that for all $\alpha, \beta \in \ob{L}$ with
             $\alpha \lcover_{\ob{L}} \beta$, we have $\mygamma_{\ob{L}} (\alpha, \beta) = 1$, then
             $\ab{A}$ is nilpotent.
       \item \label{it:tf4} If $\Con (\ab{A})$ has a $(0,1)$-sublattice $\ob{L}$ of finite height
             such that $\ob{B}_2$ is not a homomorphic image of $\ob{L}$, then $\ab{A}$ 
             is solvable.
   \end{enumerate}
\end{thm}
\emph{Proof:}
      Let $\ob{K}$ be the lattice
      $\algop{\Con(\ab{A})}{\join, \meet}$. 
      Since $\ab{A}$ lies in a congruence modular variety,
      Proposition~\ref{pro:coniscl} tells that 
        $\ab{K} := \algop{\Con (\ab{A})}{\join, \meet, [.,.]_{\ab{A}}}$ is a commutator lattice.

     \eqref{it:tf1} From Theorem~\ref{thm:forceabel}, we obtain that
       $\ob{K}$ forces abelian type, and therefore $[1_A, 1_A]_{\ab{A}} = 0_A$.
     
      \eqref{it:tf3} Let $[\alpha_1, \ldots, \alpha_n]_{\ab{A}}$ denote the 
            $n$-ary commutator of $\alpha_1,\ldots,\alpha_n \in \Con (\ab{A})$.
            A.\ Moorhead \cite{Mo:HCTF} proved that these higher commutator operations
            satisfy (among others) the conditions (HC1), (HC3), and (HC7) from
            \cite[p.~860]{AM:SOCO}. Let $c_{\ob{L}} : \ob{K} \to \ob{L}$ be the
            operation defined before Proposition~\ref{pro:ck}.
            Now for every $n \in \N$, we define an operation $f^{\ob{L}}_n : \ob{L}^n
            \to \ob{L}$ by $f^{\ob{L}}_n (x_1,\ldots, x_n) := c_{\ob{L}} ([x_1,\ldots, x_n]_{\ab{A}})$
            for all $x_1, \ldots, x_n \in \ob{L}$. Using Proposition~\ref{pro:ck},
            it is easy to verify that
             the sequence $(f^{\ob{L}}_n)_{n \in \N}$ satisfies the
            conditions (HC1), (HC3), and (HC7) of \cite{AM:SOCO}.
            Now the proof of \cite[Lemma~3.3]{AM:SOCO} yields an $n \in \N$
            with $f^{\ob{L}}_n (\underbrace{1_A,\ldots, 1_A}_n) = 0_A$, and therefore
            $[\underbrace{1_A, \ldots, 1_A}_n]_{\ab{A}} \le c_{\ob{L}} ([\underbrace{1_A,\ldots, 1_A}_n]_{\ab{A}}) =
             f^{\ob{L}}_n (\underbrace{1_A,\ldots, 1_A}_n) = 0_A$. Hence
            $\ab{A}$ is supernilpotent.

      \eqref{it:tf2} From Theorem~\ref{thm:forcenil}, we obtain that $\ob{L}$ forces nilpotent type,
        and hence by Theorem~\ref{thm:sub}, $\ob{K}$ forces nilpotent type.
        Hence  $\ab{K}$ is of nilpotent type,
        making $\ab{A}$ nilpotent.

      \eqref{it:tf4} From Theorem~\ref{thm:lattsolv}, we obtain that $\ob{L}$ forces solvable type,
        and hence by Theorem~\ref{thm:sub}, $\ob{K}$ forces solvable type.
             Thus $\ab{K}$ is of solvable type,
        making $\ab{A}$ solvable. \qed

The next sections search for partial converses of these results.
\subsection{Nonsolvable and nonnilpotent expansions}

We let $\Comp (\ab{A})$ be the clone of congruence preserving
functions of $\ab{A}$, and we define $\cc{A}$ as the algebra
$\algop{A}{\Comp (\ab{A})}$. Hence $\cc{A}$ is the largest expansion of $\ab{A}$
that has the same congruence relations as $\ab{A}$.
\begin{lem} \label{lem:3}
  Let $\ab{A}$ be an algebra in a congruence modular variety, and let $\ob{L}$ be
  its congruence lattice. We assume that $\ob{L}$ is bialgebraic.
  Let $\alpha \in J(\ob{L})$. Then $[\alpha, \alpha]_{\cc{A}} = \alpha$ if and only
  if $\alpha$ is lonesome.
\end{lem}
\emph{Proof:} For the ``only if''-direction, we assume that $[\alpha, \alpha]_{\cc{A}} = \alpha$.
Then by Lemma~\ref{lem:3lonesome},
$\alpha$ is lonesome.
For the ``if''-direction, we assume that $\alpha$ is lonesome.
Let $\eta \in M (\ob{L})$ be such that $\eta \ge \alpha^-$, $\eta \not\ge \alpha$. Then
$\Int{\alpha^-}{\alpha} \nearrow \Int{\eta}{\eta^+}$. By Proposition~\ref{pro:lonesomeJ},
$\eta$ is lonesome, and we have $\mydelta (\alpha^-, \alpha) = \alpha$ and
$\mygamma (\alpha^-, \alpha) = \eta$.
   We choose $(a,b) \in \alpha \setminus \alpha^-$ and define
   a binary function $f$ by $f (x,y) = b$ if $(x, b) \in \eta$ 
   and $(y, b) \in \eta$, and $f(x,y) = a$ else.
   By Proposition~\ref{pro:splits}, $(\eta, \alpha)$ is a splitting pair
   of the lattice $\Con (\ab{A})$, and
   the function $f$ is constant on $\eta$-classes and maps into one $\alpha$-class.
   From this we conclude that $f$ is congruence preserving (an argument is
   given, e.g., in \cite[Proposition~3.1]{ALM:FGOC}).
   Thus $f$ is a fundamental
   operation of $\cc{A}$.
   We have $f(a,a) = a = f(a,b)$, and therefore
   $\congmod{f(b, a)}{f(b, b)}{[\alpha, \alpha]_{\cc{A}}}$. Hence
   $(a, b) \in [\alpha, \alpha]_{\cc{A}}$, and therefore
   $[\alpha, \alpha]_{\cc{A}} \not\le \alpha^-$. Thus $[\alpha, \alpha]_{\cc{A}} = \alpha$. \qed

\begin{thm} \label{thm:ccAsolvable}
            Let $\ab{A}$ be an algebra in a congruence modular
            variety. We assume that $\Con (\ab{A})$ is of finite height.
            Then $\cc{A}$ is solvable if and only if $\ob{B}_2$ is not 
             a homomorphic image of the lattice
            $\Con (\ab{A})$.
\end{thm}
\emph{Proof:}
   For the ``if''-direction, we use 
   Theorem~\ref{thm:lattsolv} and obtain that $(\Con (\ab{A}), \join, \meet, [.,.]_{\ab{A}})$
    is a commutator lattice of solvable type, and therefore $\ab{A}$ is solvable.
    For the ``only if''-direction, we assume that
    $\ob{B}_2$ is a homomorphic image of $\Con (\ab{A})$.
    Let $\alpha_1, \beta_1 \in \ob{L}$ be such that $\alpha_1 \lcover \beta_1$,
    $h(\alpha_1) = 0$, and
    $h (\beta_1) = 1$. Take $\eta \in M(\ob{L})$ be such that
    $\Int{\alpha_1}{\beta_1} \nearrow \Int{\eta}{\eta^+}$. Then
    $h(\eta) = 0$ and $h(\eta^+) = 1$. Now by Proposition~\ref{pro:charlonesome},
    $\eta$ is a lonesome meet irreducible element of $\ob{L}$. Let
    $\alpha \in J(\ob{L})$ be such that $\Int {\eta}{\eta^+} \searrow \Int{\alpha^-}{\alpha}$.
    Then Proposition~\ref{pro:lonesomeJ} yields that $\alpha$ is a lonesome join irreducible
    element of $\ob{L}$.
     Now Lemma~\ref{lem:3} yields
   $[\alpha, \alpha]_{\cc{A}} = \alpha$, and hence
    by Lemma~\ref{lem:solnil}, $\cc{A}$ is not solvable. \qed

    For characterizing congruence lattices that force nilpotency,
    we restrict ourselves to finite expanded groups.
    For this characterization, we will need to construct congruence preserving
    functions that destroy nilpotency, similar to the functions
    destroying solvability
    produced in the proof of Lemma~\ref{lem:3}.
    The construction relies on certain unary congruence preserving
    functions provided by \cite{Ai:TNOC}.
    We will isolate the arguments that are restricted to expanded groups
    in the next Lemma.
   \begin{lem} \label{lem:eps}
   Let $\ab{V}$ be a finite expanded group,
   let $\ob{L}$ be its congruence lattice, and let
   $\alpha \in J(\ob{L})$.
   Then $(\alpha^- : \alpha)_{\cc{V}} \le \mygamma (\alpha^-, \alpha)$. 
  \end{lem}
    \emph{Proof:}
   We first consider the case that $\alpha$ is a lonesome
   join irreducible element. Then from Lemma~\ref{lem:3}, we
   obtain $[\alpha, \alpha]_{\cc{V}} = \alpha$. Hence
   $[\alpha, \alpha]_{\cc{V}} \not\le \alpha^-$, and therefore
   $\alpha \not\le (\alpha^- : \alpha)_{\cc{V}}$.
   Since $\alpha$ is lonesome, we have $\mydelta(\alpha⁻, \alpha) =
   \alpha$. The splitting property from Proposition~\ref{pro:splits}
    now yields $(\alpha^- : \alpha)_{\cc{V}} \le \mygamma (\alpha^-, \alpha)$.

   Let us now consider the case that $\alpha$ is not a lonesome join
   irreducible element.
 Let $A:= 0/\alpha$, $A^- := 0/\alpha^-$, 
    $C := 0/ \mygamma (\alpha^-, \alpha)$, 
    $D := 0/ \mydelta (\alpha^-, \alpha)$,
    $E := 0/ (\alpha^- : \alpha)_{\cc{V}}$.
 Our goal is to show $E \subseteq C$. To this end, we fix $z \in E$.
    We first show 
   \begin{equation} \label{eq:aga}
     \alpha \le  \mygamma (\alpha^-, \alpha).
   \end{equation}
   Since $\alpha$ is not lonesome, we apply Proposition~\ref{pro:lonesomeM}
   to the dual of $\ob{L}$ and obtain $\beta \in \Con (\ab{V})$ such that
   $\alpha$ and $\beta$ are not comparable and $\Int{\alpha^-}{\alpha} \projective
   \Int{\beta^-}{\beta}$.
   Then $\mydelta (\alpha^-, \alpha) \ge \alpha \join \beta > \alpha$,
   and therefore $\alpha \not\ge  \mydelta(\alpha^-, \alpha)$.
    Proposition~\ref{pro:splits} now yields \eqref{eq:aga}, and thus $A \subseteq C$.
    Next, we use Proposition~4.3(2)$\Rightarrow$(1) and Theorem~5.1 from \cite{Ai:TNOC}
   to obtain a unary congruence preserving function $e$ of $\ab{V}$ 
   with $e(0) = 0$, $e(A) \not\subseteq A^-$ and $e(V) \subseteq D$.
   From $e$, we define a function $f : V \times V \to V$ by 
   $f (x,y) := e (z - x + y)$ if $z - x + y \in C$ and $f(x,y) := 0$ otherwise.
   The range of $f$ is contained in $D$, and the restriction of $f$ to
   each $\mygamma (\alpha^-, \alpha)$-class, i.e., to each set of
   the form $(x_1, y_1) + C \times C$, is the restriction of a congruence
   preserving function of $\ab{V}$. By Proposition~\ref{pro:splits},
   $(\mygamma (\alpha^-, \alpha), \mydelta (\alpha^-, \alpha))$ is a splitting
   pair of $\Con (\ab{V})$, and thus from \cite[Proposition~3.1]{ALM:FGOC} we
   see that $f$ is a congruence preserving function of $\ab{V}$.
   We choose $a \in A$ such that $e(a) \not\in A^-$.
   Seeking a contradiction, we suppose that $z \not\in C$.
   We have $f(0, 0) = 0$ because $z \not\in C$ and  $f(0,a) = 0$ because $a \in C$
   and thus $z + a \not\in C$.
   Hence $(f(z, 0), f(z, a)) \in [ (\alpha^- : \alpha)_{\cc{V}}, \alpha ]_{\cc{V}}$.
   Now $f(z,0) = e(0) = 0$ because $z - z + 0 \in C$
   and $f(z, a) = f( z - z + a) = e(a)$ because $a \in C$.
   Thus $(0, e(a)) \in \alpha^-$, and therefore $e(a) \in A^-$, contradicting the choice of $a$.
   This contradiction establishes $z \in C$, which concludes the proof of $E \subseteq C$. \qed

  \begin{thm} \label{thm:cent}
     Let $\ab{A}$ be a finite expanded group, and let
     $\alpha, \beta \in \Con (\ab{A})$ be such that
     $\alpha \lcover \beta$. Then the centralizer 
     $(\alpha : \beta)_{\cc{A}}$ of $\beta$ over $\alpha$
     in $\cc{A}$ is $\mygamma (\alpha, \beta)$.
  \end{thm}
  Lemma~\ref{lem:lb} yields $\mygamma (\alpha, \beta) \le (\alpha : \beta)_{\cc{A}}$.
   Let $\alpha_1$ be minimal in $\Con (\ab{A})$ with
   $\alpha_1 \le \beta$, $\alpha_1 \not\le \alpha$.
  Then $\alpha_1$ is join irreducible, and 
   $\Int{\alpha_1^-}{\alpha_1} \nearrow \Int{\alpha}{\beta}$.
  Hence $\mygamma (\alpha, \beta) = \mygamma (\alpha_1^-, \alpha_1)$.
  From Lemma~\ref{lem:eps}, we obtain $\mygamma (\alpha_1^-, \alpha_1) \ge
  (\alpha_1^- : \alpha_1)_{\cc{A}}$, and the last expression is equal to
  $(\alpha : \beta)_{\cc{A}}$ by Lemma~\ref{lem:centproj}. This establishes the
   other inclusion. \qed

   \begin{cor} \label{cor:nil}
     Let $\ab{A}$ be a finite expanded group. Then
     $\cc{A}$ is nilpotent if and only if for all
     congruences $\alpha, \beta \in \Con (\ab{A})$ with
     $\alpha \le \beta$, we have $\mygamma(\alpha, \beta) = 1_A$.
  \end{cor}
  \emph{Proof:} For the ``if''-direction, we assume that
   for all $\alpha \lcover \beta$, we have $\mygamma(\alpha, \beta) = 1_A$.
   Then Lemma~\ref{lem:lb} implies that $(\alpha : \beta)_{\cc{A}} = 1_A$, and therefore
   $\cc{A}$ is nilpotent by Lemma~\ref{lem:solnil}. For the ``only if''-direction, we
   assume that $\cc{A}$ is nilpotent and fix $\alpha \lcover \beta \in \Con (\ab{A})$. 
   By Lemma~\ref{lem:solnil}, we then have
   $(\alpha : \beta)_{\cc{A}} = 1_A$,  and 
   Theorem~\ref{thm:cent} yields $\mygamma (\alpha, \beta) = 1_A$. \qed

We now turn to supernilpotency.
\begin{thm} \label{thm:acsnp}
   Let $\ab{A}$ be an algebra in a congruence modular
   variety, and let
   $\ob{L}$ be its congruence lattice. We assume that $\ob{L}$ is finite.
   Then $\cc{A}$ is supernilpotent if and only if $\ob{L}$ 
   does not split.
\end{thm}
\emph{Proof:}
Assume that $\ob{L}$ does not split. From \cite{Mo:HCTF},
we obtain that the higher commutator operations of $\ab{A}$
satisfy (HC1), (HC3) and (HC7) from \cite[p.\ 860]{AM:SOCO}.
  Now from the proof of
  \cite[Lemma~3.3]{AM:SOCO}, we obtain that
  $\cc{A}$ is supernilpotent. Conversely,
  assume that $(\delta, \epsi)$ is a splitting pair of 
  $\ob{L}$. Let $n \in \N$, and let $(a,b) \in \epsi$ with $a \neq b$.
  We define an $n$-ary operation
  by $f (x_1,\ldots, x_n) := a$ if at least  one of the  $x_i$ lies in 
  $a/\delta$, and $f(x_1,\ldots, x_n) := b$ else.
  Since $(\delta, \epsi)$ splits $\Con (\ab{A})$, $f$ is congruence preserving.
  Now let $y \in A \setminus (a/\delta)$. Then
  $f(y,\ldots, y) = b$.
  We use the definition of higher commutators from \cite{Bu:OTNO} (cf. \cite{AM:SAOH})
  to show that
  $[\underbrace{1,\ldots, 1}_{n \text{ times}}]_{\ab{A}} \neq 0$.
  To this end, we observe that for all $\vb{x} \in \{a,y\}^n \setminus \{(y,\ldots,y)\}$,
  we have $f(\vb{x}) = a$. Hence if $[\underbrace{1,\ldots, 1}_{n \text{ times}}]_{\ab{A}} = 0$,
  $f(y,\ldots,y,a) = f(y,\ldots, y,y)$, which means $a = b$, contradicting the choice of $a$
  and $b$.
  \qed

\section{Proofs for the Theorems from Section~\ref{sec:intro}} \label{sec:proofs}

\emph{Proof of Theorem~\ref{thm:solv}:} 
    Item~\eqref{it:s2} just spells out the definition of forcing solvability, and
    hence it is equivalent to \eqref{it:s1}.

    \eqref{it:s2}$\Rightarrow$\eqref{it:s4}:
    Let $\ab{A}$ be an algebra generating a congruence modular 
    variety with $\ob{L} \cong \Con (\ab{A})$.
  Since $\cc{A}$ can be seen as
     an expansion of
    $\ab{A}$, it generates a congruence modular variety, and we have
     $\Con (\cc{A}) \cong \ob{L}$. Thus by the assumptions, $\cc{A}$ is
     solvable.
     Now  Theorem~\ref{thm:ccAsolvable} yields that $\ob{B}_2$ is not a homomorphic image
     of $\ob{L}$.
 
    \eqref{it:s4}$\Rightarrow$\eqref{it:s2}: Let $\ab{B}$ be an algebra in
    a congruence modular variety with $\Con (\ab{B}) \cong \ob{L}$. 
    Then from Theorem~\ref{thm:lattforce}, we obtain that $\ab{B}$ is solvable. \qed

 \emph{Proof of Theorem~\ref{thm:nil}:}
    The items~\eqref{it:n1} and~\eqref{it:n2} are equivalent by the definition of forcing nilpotency.
    If \eqref{it:n2} holds and $\ab{A}$ is a finite expanded group with
    $\Con (\ab{A}) \cong \ob{L}$, then we also have $\Con (\cc{A}) \cong \ob{L}$.
    Hence from Corollary~\ref{cor:nil}, we obtain~\eqref{it:n3}.
   If \eqref{it:n3} holds, then for every finite expanded group with $\Con (\ab{B}) \cong \ob{L}$,
   Theorem~\ref{thm:lattforce} yields that $\ab{B}$ is nilpotent. \qed

 \emph{Proof of Theorem~\ref{thm:snp}:}
  The equivalence of items~\eqref{it:su1} and \eqref{it:su2} is immediate.
   \eqref{it:su2}$\Rightarrow$\eqref{it:su3}: Let $\ab{A}$ be an algebra in a congruence
   modular variety with $\Con (\ab{A}) \cong \ob{L}$. By the assumption~\eqref{it:su2}, $\cc{A}$ 
   is supernilpotent, and thus by Theorem~\ref{thm:acsnp}, $\ob{L}$ does not split. 
   \eqref{it:su3}$\Rightarrow$\eqref{it:su2}: Theorem~\ref{thm:lattforce}\eqref{it:tf3}. \qed

\section{Open Problems}
We conclude with two questions concerning congruence lattices that
make algebras abelian.
\begin{prb} 
    Characterize those modular lattices of finite height that force abelian type.
\end{prb}
 Theorem~\ref{thm:forceabel} provides one source of such lattices.
On the algebra side, a corresponding question is to describe those
lattices that force abelianity in $D$ among the lattices of finite height in $L(D)$:
\begin{prb}
   Among all lattices of finite height that are congruences lattices of some
   algebra in a congruence modular variety, characterize those $\ob{L}$ such
   that every algebra in a congruence modular variety with congruence lattice
   isomorphic to $\ob{L}$ is abelian.
\end{prb}

\section*{Acknowledgments}

The author thanks P.\ Idziak, K.\ Kearnes, and T.\ Vetterlein for valuable discussions.

\def\cprime{$'$}
\providecommand{\bysame}{\leavevmode\hbox to3em{\hrulefill}\thinspace}
\providecommand{\MR}{\relax\ifhmode\unskip\space\fi MR }
\providecommand{\MRhref}[2]{%
  \href{http://www.ams.org/mathscinet-getitem?mr=#1}{#2}
}
\providecommand{\href}[2]{#2}


\begin{thebibliography}{MMT87}

\bibitem[Aic00]{Ai:OHAH}
E.~Aichinger, \emph{On {H}agemann's and {H}errmann's characterization of
  strictly affine complete algebras}, Algebra Universalis \textbf{44} (2000),
  105--121.

\bibitem[Aic06a]{Ai:TNOC}
\bysame, \emph{The near-ring of congruence preserving functions on an expanded
  group}, Journal of Pure And Applied Algebra \textbf{205} (2006), 74--93.

\bibitem[Aic06b]{Ai:TPFO2}
\bysame, \emph{The polynomial functions of certain algebras that are simple
  modulo their center}, Contributions to general algebra. 17, Heyn, Klagenfurt,
  2006, pp.~9--24.

\bibitem[Aic14]{Ai:OTDD}
\bysame, \emph{On the {D}irect {D}ecomposition of {N}ilpotent {E}xpanded
  {G}roups}, Comm. Algebra \textbf{42} (2014), no.~6, 2651--2662.

\bibitem[ALM16]{ALM:FGOC}
E.~Aichinger, M.~Lazi{\'c}, and N.~Mudrinski, \emph{Finite generation of
  congruence preserving functions}, Monatsh. Math. \textbf{181} (2016), no.~1,
  35--62.

\bibitem[AM10]{AM:SAOH}
E.~Aichinger and N.~Mudrinski, \emph{Some applications of higher commutators in
  {M}al'cev algebras}, Algebra Universalis \textbf{63} (2010), no.~4, 367--403.

\bibitem[AM13]{AM:SOCO}
\bysame, \emph{Sequences of commutator operations}, Order \textbf{30} (2013),
  no.~3, 859--867.

\bibitem[Ava58]{Av:DSOP}
S.~P. Avann, \emph{Dual symmetry of projective sets in a finite modular
  lattice}, Trans. Amer. Math. Soc. \textbf{89} (1958), 541--558.

\bibitem[BB87]{BB:FSON}
J.~Berman and W.~J. Blok, \emph{Free spectra of nilpotent varieties}, Algebra
  Universalis \textbf{24} (1987), no.~3, 279--282.

\bibitem[BS81]{BS:ACIU}
S.~Burris and H.~P. Sankappanavar, \emph{A course in universal algebra},
  Springer New York Heidelberg Berlin, 1981.

\bibitem[Bul01]{Bu:OTNO}
A.~Bulatov, \emph{On the number of finite {M}al'tsev algebras}, Contributions
  to general algebra, 13 (Velk{\'e} Karlovice, 1999/Dresden, 2000), Heyn,
  Klagenfurt, 2001, pp.~41--54.

\bibitem[Cze08]{Cz:AOTC}
J.~Czelakowski, \emph{Additivity of the commutator and residuation}, Rep. Math.
  Logic (2008), no.~43, 109--132.

\bibitem[Cze15]{Cz:TEDC}
\bysame, \emph{The equationally-defined commutator}, Birkh\"auser/Springer,
  Cham, 2015, A study in equational logic and algebra.

\bibitem[Day69]{Da:ACOM}
A.~Day, \emph{A characterization of modularity for congruence lattices of
  algebras.}, Canad. Math. Bull. \textbf{12} (1969), 167--173.

\bibitem[FM87]{FM:CTFC}
R.~Freese and R.~N. McKenzie, \emph{Commutator theory for congruence modular
  varieties}, London Math. Soc. Lecture Note Ser., vol. 125, Cambridge
  University Press, 1987.

\bibitem[Gr{\"a}98]{Gr:GLTS}
G.~Gr{\"a}tzer, \emph{General lattice theory}, second ed., Birkh\"auser Verlag,
  Basel, 1998, New appendices by the author with B. A. Davey, R. Freese, B.
  Ganter, M. Greferath, P. Jipsen, H. A. Priestley, H. Rose, E. T. Schmidt, S.
  E. Schmidt, F. Wehrung and R. Wille.

\bibitem[Gum83]{Gu:GMIC}
H.~P. Gumm, \emph{Geometrical methods in congruence modular algebras}, vol.~45,
  Mem. Amer. Math. Soc., no. 286, American Mathematical Society, 1983.

\bibitem[Her79]{He:AAIC}
C.~Herrmann, \emph{Affine algebras in congruence modular varieties}, Acta Sci.
  Math. (Szeged) \textbf{41} (1979), no.~1-2, 119--125.

\bibitem[HH79]{HH:ACIM}
J.~Hagemann and C.~Herrmann, \emph{A concrete ideal multiplication for
  algebraic systems and its relations to congruence distributivity}, Arch.
  Math. (Basel) \textbf{32} (1979), 234--245.

\bibitem[Hig56]{Hi:GWMO}
P.~J. Higgins, \emph{Groups with multiple operators}, Proc. London Math. Soc.
  (3) \textbf{6} (1956), 366--416.

\bibitem[HM88]{HM:TSOF}
D.~Hobby and R.~McKenzie, \emph{The structure of finite algebras}, Contemporary
  mathematics, vol.~76, American Mathematical Society, 1988.

\bibitem[IS01]{IS:PRA}
P.~M. Idziak and K.~S{\l}omczy{\'n}ska, \emph{Polynomially rich algebras}, J.
  Pure Appl. Algebra \textbf{156} (2001), no.~1, 33--68.

\bibitem[Kea99]{Ke:CMVW}
K.~A. Kearnes, \emph{Congruence modular varieties with small free spectra},
  Algebra Universalis \textbf{42} (1999), no.~3, 165--181.

\bibitem[MMT87]{MMT:ALVV}
R.~N. McKenzie, G.~F. McNulty, and W.~F. Taylor, \emph{Algebras, lattices,
  varieties, volume {I}}, Wadsworth \& Brooks/Cole Advanced Books \& Software,
  Monterey, California, 1987.

\bibitem[Moo16]{Mo:HCTF}
A.~Moorhead, \emph{Higher commutator theory for congruence modular varieties},
  Talk at the conference \emph{Algebra and Algorithms} (Structure and
  Complexity Theory; A workshop on constraint satisfaction, structure theory
  and computation in algebra), University of Colorado, Boulder, May 19-22,
  2016, 2016, slides available at {\tt
  http://math.colorado.edu/algebra2016/program.html}.

\bibitem[Sco97]{Sc:TSOO}
S.~D. Scott, \emph{The structure of {$\Omega$}-groups}, Nearrings, nearfields
  and $K$-loops (Hamburg, 1995), Kluwer Acad. Publ., Dordrecht, 1997,
  pp.~47--137.

\bibitem[Smi76]{Sm:MV}
J.~D.~H. Smith, \emph{Mal'cev varieties}, Lecture {N}otes in Math., vol. 554,
  Springer Verlag Berlin, 1976.

\end{thebibliography}
\end{document}